\documentclass[11pt]{amsart}
\usepackage{amssymb,epsfig}
\usepackage{comment}
\usepackage[dvips]{color}

\newtheorem{defi}{Definition}[section]
\newtheorem{prop}[defi]{Proposition}
\newtheorem{theo}[defi]{Theorem}
\newtheorem{conj}[defi]{Conjecture}
\newtheorem{lemm}[defi]{Lemma}
\newtheorem{coro}[defi]{Corollary}
\newtheorem{rema}[defi]{Remark}
\newtheorem{exem}[defi]{Example}
\newtheorem{exems}[defi]{Examples}

\newcommand{\bdefi}{\begin{defi}}
\newcommand{\edefi}{\end{defi}}
\newcommand{\bexems}{\begin{exems}}
\newcommand{\eexems}{\end{exems}}
\newcommand{\bprop}{\begin{prop}}
\newcommand{\eprop}{\end{prop}}
\newcommand{\btheo}{\begin{theo}}
\newcommand{\etheo}{\end{theo}}
\newcommand{\blemm}{\begin{lemm}}
\newcommand{\brema}{\begin{rema}}
\newcommand{\erema}{\end{rema}}
\newcommand{\bexer}{\begin{exem}}
\newcommand{\eexer}{\end{exem}}
\newcommand{\bconj}{\begin{conj}}
\newcommand{\econj}{\end{conj}}
\newcommand{\elemm}{\end{lemm}}
\newcommand{\bcoro}{\begin{coro}}
\newcommand{\ecoro}{\end{coro}}
\newcommand{\dem}{\noindent{\bf Proof. }}
\newcommand{\rem}{\noindent{\bf Remark. }}

\newcommand{\B}{{\cal B}}

\renewcommand{\H}{{\cal H}}
\renewcommand{\O}{{\cal O}}
\newcommand{\C}{{\cal C}}

\newcommand{\maths}[1]{{\mathbb #1}}  

\newcommand{\RR}{\maths{R}}
\newcommand{\NN}{\maths{N}}
\newcommand{\CC}{\maths{C}}
\newcommand{\QQ}{\maths{Q}}

\renewcommand{\SS}{\maths{S}}

\newcommand{\HH}{\maths{H}}

\newcommand{\ZZ}{\maths{Z}}

\newcommand{\ra}{\rightarrow}

\newcommand{\ov}[1]{{\overline #1}} 
\newcommand{\wt}[1]{{\widetilde{#1}}}

\newcommand{\ga}{\gamma}
\newcommand{\Ga}{\Gamma}

\newcommand{\cqfd}{\hfill$\Box$}

\newcommand{\hdr}{{\HH\,}^2_\RR}

\newcommand{\hnr}{{\HH\,}^n_\RR}

\newcounter{fig}

\newtheorem{lemma}[defi]{Lemma}
\newtheorem{proposition}[defi]{Proposition}
\newtheorem{theorem}[defi]{Theorem}
\newtheorem{corollary}[defi]{Corollary}

\theoremstyle{definition}
\newtheorem{remark}[defi]{Remark}

\newcommand{\hdc}{{\HH\,}^2_\CC}
\newcommand{\cal}{\mathcal}
\newcommand{\height}{\operatorname{ht}}

\begin{document}
\title{Unclouding the sky of negatively curved manifolds} 
\author{Jouni Parkkonen \and Fr\'ed\'eric Paulin}


\begin{abstract}
  Let $M$ be a complete simply connected Riemannian manifold, with
  sectional curvature $K\leq -1$. Under some assumptions on the
  geometry of $\partial M$, which are satisfied for instance if $M$ is
  a symmetric space, or has dimension $2$, we prove that given any
  family of horoballs in $M$, and any point $x_0$ outside these
  horoballs, it is possible to shrink uniformly, by a finite amount
  depending only on $M$, these horoballs so that some geodesic ray
  starting from $x_0$ avoids the shrunk horoballs. As an application,
  we give a uniform upper bound on the infimum of the heights of the
  closed geodesics in the finite volume quotients of $M$.
\end{abstract}

\maketitle

\section{Introduction}

Let $M$ be a complete simply connected Riemannian manifold, with
negative sectional curvature $K\leq -1$, and $\partial M$ be its space
at infinity, i.e.~the set of asymptotic classes of geodesic rays in
$M$. A horoball $H\!B$ in $M$ can be defined as a proper non-empty
subset of $M$ which is the limit (for the Hausdorff distance on
compact subsets) of a sequence of closed balls of radius $r_n$
converging to $\infty$.  For every $t\geq 0$, define the shrunk
horoball $H\!B(t)$ to be the horoball which is the limit of the
sequence of balls with the same center and radius $r_n -t$. Given a
collection $(H\!B_n)_{n\in\NN}$ of horoballs with pairwise disjoint
interiors, and a point $x_0$ in $M$ or $\partial M$, it may happen
that every geodesic ray or line starting from $x_0$ meets at least one
of the horoballs $H\!B_n$. We are interested in proving that there
exists a finite and explicit lower bound on $t$ such that at least one
geodesic ray or line starting from $x_0$ avoids the uniformly shrunk
horoballs $H\!B_n(t)$. The main point of this paper (besides the
existence) is to get universal (and possibly as small as possible)
such lower bounds. Thinking of $\partial M$ as the sky of $M$, and of
horoballs as clouds, we are interested in shrinking the clouds to be
able to see the blue sky.

\btheo \label{theointro}
Let $X$ be a proper 
CAT$(-1)$ metric space, and
$t_{\rm min}>0$. Assume that one of the following conditions holds:
\begin{enumerate}
\item $X$ is the real hyperbolic $n$-space with $n\geq 2$, and
  $t_{\rm min}=-\log(4\sqrt{2}-5)$;
\item $X$ is a complete simply connected Riemannian manifold of
  dimension $2$, with 
  curvature $-a^2\leq K\leq -1$, and
  $e^{-t_{\rm min}}=
  2^{2/a}\left(\sqrt{1+2^{1-1/a}}-1-2^{-1-1/a}\right)$, where $1\leq a\leq
  2$;
\item $X$ is a locally finite tree (without vertices of degree $1$ or
  $2$, with edge lengths $1$) and $t_{\rm min}=1$; 
\item $X$ is a negatively curved symmetric space, and $t_{\rm min}$ is
  some constant depending only on $X$.
\end{enumerate}
Let $(H\!B_n)_{n\in\NN}$ be a sequence of horoballs with pairwise
disjoint interiors.  If $t> t_{\rm min}+2$, then for every $x_0$ in
$X-\bigcup_n H\!B_n$, there exists at least one geodesic ray starting
at $x_0$, which avoids $H\!B_n(t)$ for every $n$ in $\NN$.

Furthermore, if $t> t_{\rm min}$, then there exists at
least one geodesic line starting at the point at infinity $x_0$ of
$H\!B_0$ which avoids $H\!B_n(t)$ for every $n$ in $\NN-\{0\}$.
\etheo

When $x_0$ is at infinity, the first three cases are sharp results.
In Sections \ref{sec:spaboundsha} and \ref{sec:bi-infinite}, we give
sufficient conditions on the geometry of the space at infinity of $X$
for the results to be valid.  The complete result is more general than
the one above. In particular, the second result of Theorem
\ref{theointro} holds for $X$ a homogeneous negatively curved
Riemannian manifold of Carnot type (see Section
\ref{sec:spaboundsha}), when $x_0$ is identified with the only point
fixed by a simply transitive group of isometries of $X$.
Note that the result is not true for every proper geodesic CAT$(-1)$ metric
space, as examples of trees with unbounded edge lengths show.
%

\medskip Consider now a finite volume complete negatively curved
Riemannian manifold $V$, and $e$ an end of $V$. Let $\beta_{e}$ be the
Busemann function on $V$ with respect to $e$, normalized to be $0$ on
the boundary of the maximal Margulis neighbourhood of $e$ (see section
\ref{sec:bi-infinite}).  Define the height of any compact subset $A$
of $V$ as the maximum value of $\beta_{e}$ on $A$.  As an application
of Theorem \ref{theointro}, we give a uniform upper bound, which
depends only on the universal cover of $V$, on the infimum $h_e(V)$ of
the heights of closed geodesics on $V$.

\bcoro \label{corointro}
Let $t_{\rm min}>0$ and let $V$ be a finite volume complete
Riemannian manifold, with negative sectional curvature $K\leq -1$.
If $t_{\rm min}$ and the universal cover of $V$  
satisfy the condition (1), (2) or (4) of the
above theorem,   then $h_e(V)\leq t_{\rm min}$.
\ecoro

The results of this paper are related to the problem of finding
geodesic rays avoiding obstacles in  negatively curved
manifolds (and constructing proper invariant subsets for the geodesic
flow, see for instance the work of Burns and Pollicott \cite{BP} and
Buyalo, Schroeder and Walz \cite{BSW}) or to
the shrinking target problem of Hill and Velani \cite{HV}.  They are
also related to the study of bounded geodesics in complete finite
volume negatively curved manifolds, see for instance \cite{Dan,Sch}
and the references therein. In particular, V.~Schroeder \cite{Sch}
proved the existence of a geodesic line through any given point having
bounded height (though not with an effective bound, and in dimension at
least $3$, but without any other assumption besides having sectional
curvature at most $-1$).
\bigskip

Let us describe briefly the contents of this paper.  In Section
\ref{sec:geomhoro}, we collect geometric results on the shadows (in
the sense of Sullivan
) of balls and horoballs in CAT($-1$) spaces in
terms of the natural distance (due to Hamenst\"adt \cite{Ham1}, see
\cite[Appendix]{HP1}) in the punctured boundary.

In Section \ref{sect:abstractuncov}, we prove an abstract uncovering
result for collections of balls in metric spaces.  

\bprop \label{propintro} 
Let $Y$ be a complete locally compact metric space such that through
two points passes a geodesic line. Then for every family
$(B(x_n,r_n))_{n\in\NN}$ of balls in $Y$, with $0< r_n\leq 1$ and
$d(x_n,x_m)^2\geq 4r_nr_m$ for $n\neq m$, the scaled family of balls
$(B(x_n,sr_n))_{n\in\NN}$ no longer covers $Y$ if $s<\sqrt{5}-2$.
\eprop

The fact that there is such a universal constant is amazing. In fact,
the result is true, with different bounds on the scaling constant $s$,
for a larger class of spaces which includes the Heisenberg group with
its Carnot-Carath\'eodory metric.  See Theorem \ref{thm:abstract} for
a more precise statement.  This result (or variations on its proof)
will imply Theorem \ref{theointro}, as in Section \ref{sec:geomhoro},
we prove that the disjointness of horoballs implies the above
quadratic separation property for the packing of balls on the boundary
at infinity obtained from the shadows of the horoballs.

In Section \ref{sec:spaboundsha} we study the metric properties at
infinity of negatively curved Riemannian manifolds, whose punctured
boundaries are Carnot groups (see \cite{FS,Gro,HK} for instance). We
show that these manifolds satisfy the sufficient conditions needed for
Theorem \ref{theointro} to be applicable.

The improvement of the constants, in order to get sharp results, is
done for two-dimensional manifolds with pinched negative curvature in
Section \ref{sec:2-dim}, and extended to the real hyperbolic space by
a symmetry trick in Section \ref{sec:hnr}.

Finally, in Section \ref{sec:bi-infinite}, we turn to the case of
equivariant families of horoballs. This is the motivating problem, in
particular for families of horoballs which arise from arithmetic
constructions (see for instance Section \ref{sec:2-dim}). We prove
Corollary \ref{corointro} in this last section by showing the
existence of a geodesic which does not enter deep in the lift of a
maximal Margulis neighbourhood of the end $e$ and applying results of
\cite{HP2}.  
The uniformity of the upper bound on $h_e(V)$ is surprising. The value
$t_{\rm min}= -\log(4\sqrt{2}-5)\approx 0.4205$ of this upper bound in
the constant curvature case is even not very far from the bound in
particular cases coming from arithmetic constructions. For instance,
if $V$ is the orbifold ${\rm PSL}_2(\RR)\backslash \HH^2_\RR$, which
has only one end $e$, then $h_e(V)= \log (\sqrt{5}/2)\approx
0.1115$. This has been known for a long time, as $1/\sqrt{5}$ is the
Hurwitz constant of the classical diophantine approximation problem of
real numbers by rational numbers, but see \cite{HP3} and the
references therein, and also \cite{Ser,RWT}).

\medskip
\noindent
{\small
{\it Acknowledgements: }
We thank P.~Pansu for explaining us the dilation trick used to
compare the Hamenst\"adt and Carnot-Carath\'eodory metrics on the
boundary of negatively curved homogeneous spaces in section
\ref{sec:spaboundsha}. The first author also acknowledges the support
of the Center of Excellence "Geometric analysis and mathematical
physics" of the Academy of Finland. Part of this work was completed
during a stay of the first author at 
the Ecole Normale Sup\'erieure.
}

\section{The geometry of shadows in negatively curved manifolds}
\label{sec:geomhoro}

We refer to \cite{Bourdon,BH,GH} for the definitions and the first
properties of a CAT$(-1)$ geodesic metric space $X$ and its boundary
$\partial X$, as well as the horospheres and horoballs in X.  All
balls and horoballs are assumed to be closed unless otherwise stated.
Accordingly, in any metric space $(Y,\delta)$, we will denote by
$B(x,r)=B_\delta(x,r)$ the closed ball of center $x$ and radius $r$,
and by $S(x,r)=S_\delta(x,r)$ the sphere of center $x$ and radius $r$.
If $\xi,\eta$ are two points in $\partial X$, we will denote by
$]\xi,\eta[$ the geodesic line in $X$ between $\xi$ and $\eta$.

\bigskip 
Let $(X,d)$ be a CAT($-1$) geodesic metric space, let
$\xi_\sharp$ be a point in $\partial X$, and let $H_\sharp$ be a
horosphere in $X$ centered at $\xi_\sharp$.  The {\it Hamenst\"adt
  distance} $d_{\xi_\sharp ,H_\sharp}$ on $\partial X-\{\xi_\sharp\}$
is defined by
\begin{equation}\label{eqn:defhamdist}
d_{\xi_\sharp,H_\sharp}(a,b)=\lim_{t\ra +\infty}
e^{-\frac{1}{2}(2t-d(a_t,b_t))}
\end{equation}
for $t\mapsto c_t$ the geodesic line from $\xi_\sharp$ to some point
$c$ in $\partial X-\{\xi_\sharp\}$, with $c_0$ belonging to
$H_\sharp$. We refer to \cite[Appendix]{HP1}\cite{HP2} for more
details.

For instance, if $X$ is the upper halfplane model of the real
hyperbolic $n$-space $\hnr$, if $\xi_\sharp$ is the point at infinity
$\infty$, and if $H_\sharp$ is the horosphere at Euclidean height $1$,
then the Hamenst\"adt distance $d_{\infty,H_\sharp}$ is the Euclidean
distance on $\partial\hnr-\{\infty\}=\RR^n$.

For $H'_\sharp$ another horosphere centered at $\xi_\sharp$, denote by
$d_{alg}(H'_\sharp, H_\sharp)$ the algebraic distance, on any geodesic
line $\ell$ ending at $\xi_\sharp$ (and oriented towards
$\xi_\sharp$), between the intersection points of $\ell$ with
$H'_\sharp$ and $H_\sharp$. We then have
\begin{equation}\label{eqn:disthamvarihoro}
d_{\xi_\sharp,H'_\sharp}=e^{-d_{alg}(H'_\sharp, H_\sharp)}
\,d_{\xi_\sharp,H_\sharp}\;.
\end{equation}  
Note that if $\ga$ is an isometry of $X$, then, for every $a,b$ in
$\partial M-\{\xi_\sharp\}$,
\begin{equation}\label{eqn:invadistham}
d_{\ga \xi_\sharp,\ga H_\sharp}(\ga a,\ga b)=
d_{\xi_\sharp,H_\sharp}(a,b)\;.
\end{equation}

\medskip\noindent{\bf Remark.} The Hamenst\"adt distance has the
following scaling property.  For every $\epsilon$ in $(0,1]$, the
metric space $\epsilon X=(X,\epsilon d)$ is again a CAT$(-1)$ geodesic
metric space, with boundary $\partial(\epsilon X)$ naturally
identified with $\partial X$, and $H_\sharp$ is again a horosphere in
$\epsilon X$, centered at the point at infinity $\xi_\sharp$ of
$\epsilon X$.  It follows by definition that the Hamenst\"adt distance
in $\partial(\epsilon X)$ defined by $\xi_\sharp$ and $H_\sharp$ is
exactly ${d_{\xi_\sharp,H_\sharp}}^\epsilon$.

\brema \label{rem:hamlocequiv}
{\rm The Hamenst\"adt distances with respect to two points at
  infinity $\xi_\sharp,\xi'_\sharp$ are locally equivalent. 
More precisely, for every
  $\xi_\sharp,\xi'_\sharp,\eta_0$ in $\partial X$, for all
  horospheres $H_\sharp,H'_\sharp$ centered at
  $\xi_\sharp,\xi'_\sharp$ respectively, there exists a constant
  $\lambda>0$  and a neighborhood
  $U$ of $\eta_0$ such that for every 
  $\eta,\eta'$ in $U$, we have
  $$\lambda \;e^{-8}\;
  d_{\xi_\sharp,H_\sharp}(\eta,\eta')\leq
  d_{\xi'_\sharp,H'_\sharp}(\eta,\eta') \leq
  \lambda \;e^8\;d_{\xi_\sharp,H_\sharp}(\eta,\eta')\;.$$
}
\erema

\dem 
Assume that $\xi_\sharp,\xi'_\sharp,\eta_0,H_\sharp,H'_\sharp$
are as in the statement.  Let $x$ be a point on $]\xi_\sharp,\eta_0[$
at distance at most $1$ from $]\xi'_\sharp,\eta_0[$ and
$]\xi_\sharp,\xi'_\sharp[$, which exists by an easy comparison
argument with the hyperbolic plane. Then for every small enough
neighborhood $U$ of $\eta_0$ and for every $\eta$ in $U$, the geodesic
lines between $\xi_\sharp$ and $\eta$, and between $\xi'_\sharp$ and
$\eta$, are passing at distance at most $2$ from the point $x$.  Let
$\eta,\bar\eta$ be in $U$. Let $p,\bar p$ (resp.~$p',\bar p'$) be the
closest points to $x$ on $]\xi_\sharp,\eta[$ and
$]\xi_\sharp,\bar\eta[$ (resp.~on $]\xi'_\sharp,\eta[$ and
$]\xi'_\sharp,\bar\eta[$). We chose $U$ such that $p,\bar p,p',\bar
p'$ belong to $B(x,2)$. Up to changing the horoballs
$H_\sharp,H'_\sharp$, which will only change the value of $\lambda$ by
Equation \eqref{eqn:disthamvarihoro}, we may assume that $H_\sharp$
(resp.~$H'_\sharp$) contains the point $p$ (resp.~$p'$).  The triangle
inequality implies
$$
\begin{aligned}
& \left| \left( d(H_\sharp,\eta_t) + d(H_\sharp,\bar\eta_t) -
  d(\eta_t,\bar\eta_t)\right)-\left( d(x,\eta_t) + d(x,\bar\eta_t) -
  d(\eta_t,\bar\eta_t)\right) \right|\\ 
& \leq  d(H_\sharp,p)+d(H_\sharp,\bar p)+d(p,x)+d(\bar p,x)\leq 
  d(\bar p,x)+d(x,p)+d(p,x)+d(\bar p,x)\leq 8\;,
\end{aligned}
$$ 
and similarly
with $H'_\sharp$ replacing $H_\sharp$. Therefore 
$$e^{-8}\;\;d_{\xi_\sharp,H_\sharp}(\eta,\bar\eta)\leq
d_{\xi'_\sharp,H'_\sharp}(\eta,\bar \eta) \leq
e^8\;d_{\xi_\sharp,H_\sharp}(\eta,\bar\eta)\;.
\;\;\;\;\;\;\;\;\mbox{\cqfd}$$

If $A$ is a subset of $X$, define the {\it shadow} of $A$ seen from
$\xi_\sharp$ to be the subset $\O A=\O_{\xi_\sharp}(A)$ of $\partial
X-\{\xi_\sharp\}$ consisting of the endpoints of the geodesic lines
starting from $\xi_\sharp$ and meeting $A$. Note that if $A$ is the
(closed) horoball bounded by a horosphere $H$, then $\O H=\O A$, as
$H$ separates $X$.

It is geometrically clear that, in the upper halfspace model of the
real hyperbolic $n$-space, the shadow seen from $\infty$ of a ball or
horoball is a Euclidean ball. Our first result says that when $X$ is
assumed to be a pinched negatively curved Riemannian manifold, then
the shadow of a ball or horoball is almost a ball for the Hamenst\"adt
distance. Comparing shadows of balls with balls for some distance on
the boundary has a long history, starting with the shadow lemma of
\cite{Sul}, see also \cite{Bourdon,Rob} and \cite[Lemma 3.1]{HP3}.

\medskip
\begin{proposition} \label{prop:horosizshad} 
  Let $M$ be a complete simply connected Riemannian manifold, with
  sectional curvature $-a^2\leq K\leq -1$, let $\xi_\sharp$ be a point
  at infinity and $H_\sharp$ a horosphere centered at $\xi_\sharp$.
  For every ball $B$ of center $x$ and radius $r$, whose interior is
  disjoint from the open horoball bounded by $H_\sharp$, we have
$$
B_{d_{\xi_\sharp,H_\sharp}}(\xi, 
\frac{1}{2}e^{-d(H_\sharp,B)}(1-e^{-2r}))
\subset \O_{\xi_\sharp}(B) \subset
B_{d_{\xi_\sharp,H_\sharp}}(\xi,2^{-\frac{1}{a}}
e^{-d(H_\sharp,B)}(1-e^{-2ar})^\frac{1}{a})\;,
$$
where $\xi$ is the other endpoint of the geodesic line from
$\xi_\sharp$ through $x$.
\end{proposition}

\medskip
\noindent{\bf Remark.}  Note that this result is sharp, since the
proof shows that if $M$ has constant curvature $-1$, then the left
inclusion is an equality, and if $M$ has constant curvature $-a^2$,
then the right inclusion is an equality.

\bigskip 
\dem
To simplify notations, let $d_\infty= d_{\xi_\sharp,H_\sharp}$. For
every $\kappa\geq 1$, let us denote by $\HH^2_{-\kappa}$ the upper
halfplane model of the real hyperbolic plane with constant curvature
$-\kappa$, and $d_{-\kappa}$ its distance. Let $\ov{\xi_\sharp}$ be
the point at infinity in this model, and let $\ov{H_\sharp}$ be the
horosphere centered at $\ov{\xi_\sharp}$ at Euclidean height $1$. Let
$d_{-\kappa,\infty}$ be the Hamenst\"adt distance on
$\partial\HH^2_{-\kappa}-\{ \ov{\xi_\sharp}\}$ with respect to the
horoball $\ov{H_\sharp}$.

\bigskip 
Let us start with a preliminary remark. Let $\ov{\xi},\ov{\xi'}$ be
two distinct points in $\partial \HH^2_{-\kappa}-\{\ov{\xi_\sharp}\}$,
and let $\ov{x}$ be a point in $]\ov{\xi_\sharp},\ov{\xi}[$ below
$\ov{H_\sharp}$.  Define $\ov{p'}$ to be the orthogonal projection of
$\ov{x}$ on $]\ov{\xi_\sharp},\ov{\xi'}[$.  Define $\ov{B'}$ to be the
ball centered at $\ov{x}$ whose boundary contains $\ov{p'}$.
Denote its radius by $\overline{r}$. Assume that the interior of
$\ov{B'}$ is disjoint from the open horoball bounded by
$\ov{H_\sharp}$ (see Figure \ref{fig:inclusionombreboule} below).

First assume that $\kappa=1$, hence in particular the Hamenst\"adt
distance $d_{-1,\infty}$ is the Euclidean distance on
$\partial\HH^2_{-1}-\{ \ov{\xi_\sharp}\}=\RR$.  Denote by $u$ the
Euclidean height of $\ov{x}$ and by $R$ the Euclidean radius of the
circle containing the geodesic segment between $\ov{x}$ and $\ov{p'}$.
Then
$$u=e^{-(\ov{r}+d_{-1}(\ov{H_\sharp},\ov{B'}))}$$
and
$$R^2=u^2 +d_{-1,\infty}(\ov{\xi},\ov{\xi'})^2$$
by Pythagoras' formula in the right angled Euclidean triangle
$\ov{x}\ov{\xi}\ov{\xi'}$, and
$$R+d_{-1,\infty}(\ov{\xi},\ov{\xi'})=
e^{-d_{-1}(\ov{H_\sharp},\ov{B'})}$$
as the Euclidean height of $\ov{p'}$ is $R$ and the Euclidean center
of the hyperbolic ball $\ov{B'}$ has the same Euclidean height as
$\ov{p'}$. Eliminating $u$ and $R$, we get
$$ d_{-1,\infty}(\ov{\xi},\ov{\xi'})= 
\frac{1}{2}e^{-d_{-1}(\ov{H_\sharp},\ov{B'})}(1-e^{-2\ov{r}})\;.
$$

Assuming now that $\kappa\geq 1$, by the scaling property of the
Hamenst\"adt distance, and as $\HH^2_{-\kappa}=
\frac{1}{\sqrt{\kappa}}\HH^2_{-1}$, we have
$$
d_{-\kappa,\infty}(\ov{\xi},\ov{\xi'})= \left(
\frac{1}{2}\,e^{-\sqrt{\kappa}\,d_{-\kappa}(\ov{H_\sharp},\ov{B'})}
(1-e^{-2\sqrt{\kappa}\,\ov{r}})\right)^\frac{1}{\sqrt{\kappa}}
=2^{-\frac{1}{\sqrt{\kappa}}}\,e^{-d_{-\kappa}(\ov{H_\sharp},\ov{B'})}
(1-e^{-2\sqrt{\kappa}\,\ov{r}})^\frac{1}{\sqrt{\kappa}}\;.
$$

\bigskip 
Let us first prove the inclusion on the right in the statement of
Proposition \ref{prop:horosizshad}. Take $\xi'$ in $\O B$. Let us
prove that $d_\infty(\xi,\xi')\leq 2^{-\frac{1}{a}}
e^{-d(H_\sharp,B)}(1-e^{-2ar}) ^\frac{1}{a}$.


\begin{figure}[ht]
\input{fig_jouni20fev03new.pstex_t}
\caption{}\label{fig:inclusionombreboule}
\end{figure}

We may assume that $\xi$ and $\xi'$ are different. Let $p'$ be the
orthogonal projection of $x$ on $]\xi_\sharp,\xi'[$, that is, the
closest point to $x$ on the geodesic line $]\xi_\sharp,\xi'[$. Let
$B'$ be the ball in $M$ centered at $x$ whose boundary contains $p'$,
and hence is tangent to $]\xi_\sharp,\xi'[$ at $p'$. As
$]\xi_\sharp,\xi'[$ meets $B$ and $p'$ is the closest point, the ball
$B'$ is contained in $B$.

There exists two distinct points $\ov{\xi},\ov{\xi'}$ in $\partial
\HH^2_{-a^2}-\{\ov{\xi_\sharp}\}$, and a point $\ov{x}$ in
$]\ov{\xi_\sharp},\ov{\xi}[$ below $\ov{H_\sharp}$, such that
$$d_\infty(\xi,\xi')=d_{-a^2,\infty}(\ov{\xi},\ov{\xi'})\;\;{\rm
  and}\;\; d(x,H_\sharp)= d_{-a^2}(\ov{x},\ov{H_\sharp})\;.$$
Define $\ov{p'},\ov{B'},\ov{r}$ as in the preliminary remark. By the
comparison property, 
as the sectional curvature of $M$ is at least $-a^2$, by taking limits
of comparison triangles (see for instance \cite{GH}), 
we have
$$\ov{r}=d_{-a^2}(\ov{x},\ov{p'})\leq d(x,p')\leq r\;,$$ which implies
that
$$d_{-a^2}(\ov{H_\sharp},\ov{B'})\geq d(H_\sharp,B')\geq
d(H_\sharp,B)\;.$$ In particular, the interior of $\ov{B'}$ is
disjoint from the open horoball bounded by $\ov{H_\sharp}$. Hence, by
the preliminary remark,
$$d_\infty(\xi,\xi')=d_{-a^2,\infty}(\ov{\xi},\ov{\xi'})=
2^{-\frac{1}{a}}\,e^{-d_{-a^2}(\ov{H_\sharp},\ov{B'})}
(1-e^{-2a\ov{r}})^\frac{1}{a}\leq 2^{-\frac{1}{a}}\,e^{-d(H_\sharp,B)}
(1-e^{-2ar})^\frac{1}{a}\;.$$
This is what we wanted.

\medskip 
The proof of the other inclusion is similar (though easier, as we
compare only with the constant curvature $-1$ hyperbolic plane), and
is left to the reader.
\cqfd

\medskip
\begin{corollary} \label{coro:ballshadsize} 
Under the assumption of Proposition \ref{prop:horosizshad}, for every
horosphere $H$ centered at a point at infinity $\xi$, bounding an open
horoball disjoint from the open horoball bounded by $H_\sharp$, we
have
$$
B_{d_{\xi_\sharp,H_\sharp}}(\xi,\frac{1}{2}e^{-d(H_\sharp,H)})
\subset \O_{\xi_\sharp}H \subset
B_{d_{\xi_\sharp,H_\sharp}}(\xi,2^{-\frac{1}{a}}
e^{-d(H_\sharp,H)})\;.
$$
As previously, this result is sharp.
\end{corollary}

\dem Let $p$ be the intersection point of $H$ and $]\xi_\sharp,\xi[$.
Let $x_t$ be the point on $]\xi_\sharp,\xi[$, at distance $t> 0$ from
$p$ (and converging to $\xi$ as $t$ tends to $+\infty$). Let $B_t$ be
the ball of center $x_t$ whose boundary contains $p$. Then $B_t$
converges, for the topology of uniform convergence on compact subsets
of $M$, to the horoball $H\!B$ bounded by $H$.  Hence $\O B_t$
converges, for the Hausdorff distance on closed subsets of $\partial
M-\{\xi_\sharp\}$ to $\O H$. Note that the radius of $B_t$ tends to
$+\infty$. Hence the result follows from Proposition
\ref{prop:horosizshad} by taking limits.
\cqfd

\medskip
\begin{remark}\label{rem:ombreboulegene}{\rm 
It is well-known that some version of this sharp result holds for
general CAT$(-1)$ spaces, and here is one, where the constant $c_{\rm
max}$ is probably not optimal. We claim that for 
$c_{\rm max}=e^{{2}}$, 
with the standing assumptions in this section
on $(X,\xi_\sharp,H_\sharp)$, for every horosphere $H$ in $X$ centered
at a point at infinity $\xi$, bounding an open horoball disjoint from
the open horoball bounded by $H_\sharp$, we have
$$
B_{d_{\xi_\sharp,H_\sharp}}(\xi,\frac{1}{2}e^{-d(H_\sharp,H)})
\subset \O_{\xi_\sharp}H \subset
B_{d_{\xi_\sharp,H_\sharp}}(\xi,c_{\rm max}\, e^{-d(H_\sharp,H)})\;.
$$
\dem The left inclusion is the one in the corollary above, since its
proof only used the fact that $M$ is CAT$(-1)$. To prove the right
inclusion, let $\eta$ be the endpoint of a geodesic line starting from
$\xi_\sharp$ and meeting $H$, that we may assume to be different from
$\xi$. Let $p$ be the orthogonal projection of $\eta$ to
$]\xi_\sharp,\xi[$, which is the closest point to $\eta$ on the
geodesic line between $\xi_\sharp$ and $\xi$ in the sense of any
Busemann function corresponding to $\eta$. Note that as in Proposition
\ref{prop:horosizshad}, the point $p$ belongs to the horoball bounded
by $H$. Let $q$ (resp.~$q_\sharp$) be the point, which is the closest
to $p$, on the geodesic line between $\eta$ and $\xi$ (resp.~between
$\eta$ and $\xi_\sharp$). By an easy computation in the
upper-halfspace model, and by comparison using the CAT$(-1)$ property
of $X$, the distance between $p$ and $q$, and the one between $p$ and
$q_\sharp$, are at most $1$. Let $t\mapsto \xi_t$
(resp.~$t\mapsto\eta_t$) be the geodesic line starting from
$\xi_\sharp$, passing through $H_\sharp$ at time $t=0$, and ending in
$\xi$ (resp.~$\eta$). For every $\epsilon>0$, for $t$ big enough, as
the geodesic segment between $\xi_t$ and $\eta_t$ converges to the
geodesic line between $\xi$ and $\eta$, there exists a point on
$[\xi_t,\eta_t]$ which is at distance at most $\epsilon$ from $q$.
Hence, by applying the triangular inequality several times, we have
$$
\begin{aligned}
 & d(\xi_t,H_\sharp) + d(\eta_t,H_\sharp) - d(\xi_t,\eta_t)\\ & \geq
 d(\xi_t,p)+d(p,\xi_0)+d(\eta_t,q_\sharp)+d(q_\sharp,\eta_0)
    -d(\xi_t,q)-d(q,\eta_t)\\ 
&\geq -d(p,q)-d(q_\sharp,q)+d(\xi_0,p)+d(\eta_0,q_\sharp)\geq 
2d(\xi_0,p)-4\;.
\end{aligned}
$$
Therefore
$$d_{\xi_\sharp,H_\sharp}(\eta,\xi)\leq 
e^{-d(p,H_\sharp)+{2}}\leq
c_{\rm max}\, e^{-d(H,H_\sharp)} \;.\;\;\;\;\mbox{\cqfd}
$$
}
\end{remark}

\bigskip 
If $H$ is a horosphere in $X$, bounding an open horoball disjoint
from $H_\sharp$, motivated by Proposition \ref{prop:horosizshad} and
Remark \ref{rem:ombreboulegene}, we call
$$ 
 r_i(H)=\frac 12 e^{-d(H_\sharp,H)} 
$$
the {\it inner radius of the shadow} of $H$. When $X$ is furthermore a
complete Riemannian manifold with sectional curvature $-a^2\leq K\leq
-1$ (with $a\geq 1$), we call
$$
  r_e(H)=2^{-\frac 1a} e^{-d(H_\sharp,H)}
$$
the {\it outer radius of the shadow} of $H$. These quantities depend
on $H_\sharp$.  Note that the outer radius is bounded above by
$e^{-d(H_\sharp,H)}$.

\bigskip
We are now going to give an inequality relating the inner radii of
two horospheres bounding disjoint open horoballs to the Hamenst\"adt
distance between their points at infinity.

\medskip
Given two balls $B=B(x,r)$ and $B'=B(x',r')$ in $X$, define the {\it
algebraic distance} between $B$ and $B'$ to be
\begin{equation}
d_{\rm alg}(B,B')= d(x,x')-r -r'\;.
\label{eqn:distalgball}
\end{equation}
Note that $d_{\rm alg}(B,B')=d(B,B')\geq 0$ if $B$ and $B'$ have
disjoint interior, as $X$ is geodesic.

Given two points $x,x'$ in $X$, that do not belong to the open
horoball bounded by $H_\sharp$, define
\begin{equation}
d_{\xi_\sharp,H_\sharp}(x,x')= e^{-\frac{1}{2}(d(x,H_\sharp)
+d(x,H_\sharp)-d(x,x'))} \;.
\label{eqn:disthampoints}
\end{equation}
Note that, by the definition of the Hamenst\"adt distance, see
Equation (\ref{eqn:defhamdist}), the numbers
$d_{\xi_\sharp,H_\sharp}(x_i,x'_i)$ converge to
$d_{\xi_\sharp,H_\sharp}(\xi,\xi')$ as the points $x_i,x'_i$ tend
respectively to the points $\xi,\xi'$ in $\partial X-\{\xi_\sharp\}$
while staying on a geodesic ray that converges to $\xi,\xi'$.

\blemm\label{lem:hamdistballs} 
Let $B, B'$ be two balls in $X$ with center $x, x'$, whose interior do
not meet the open horoball bounded by $H_\sharp$. Then
$$
-2\log d_{\xi_\sharp,H_\sharp}(x,x') = d(H_\sharp,B) +
 d(H_\sharp,B')-d_{alg}(B,B')\;
$$
\elemm

\dem If $r$ is the radius of $B$ then $d(H_\sharp,B)=d(H_\sharp,x)-r$.
The result follows from the equations (\ref{eqn:distalgball}) and
(\ref{eqn:disthampoints}). 
\cqfd

\medskip
If $H, H'$ are two horospheres in $X$ centered at the distinct points
$a, b$ in $\partial X$, denote by $d_{\rm alg}(H,H')$ the distance
between the intersection points with $H$ and $H'$ of the geodesic line
between $a$ and $b$, with a positive sign if $H$ and $H'$ are
disjoint, and a negative sign if the open horoballs bounded by $H$ and
$H'$ meet.

\bcoro\label{coro:disthamhorob}
Let $H, H'$ be two horospheres in $X$ centered at distinct points
$\xi, \xi'$ in $\partial X$ respectively, with the open horoball
bounded by $H_\sharp$ disjoint from the open horoballs bounded by $H,
H'$. Then
$$
-2\log d_{\xi_\sharp,H_\sharp}(\xi,\xi') = d(H_\sharp,H) + d
 (H_\sharp,H')-d_{alg}(H,H')\;
$$
\ecoro

\dem
This follows from Lemma \ref{lem:hamdistballs} by taking
limits, as in the proof of Corollary \ref{coro:ballshadsize}. 
\cqfd

\medskip
\bcoro\label{coro:disjointhamenstadt} 
With the notations of the previous corollary, assume furthermore that
$H$ and $H'$ bound disjoint open horoballs. Then
\begin{equation}
 d_{\xi_\sharp,H_\sharp}(\xi,\xi')^2\geq 4
 r_i(H)r_i(H'). \label{eqn:disjhamineq} 
\end{equation}
This inequality is sharp, i.e.~the equality holds if and only if $H$
 and $H'$ meet in one (and only one) point (their common intersection
 point with the geodesic line between $\xi$ and $\xi'$). 
\ecoro

\dem 
By Corollary \ref{coro:disthamhorob}, since the open horoballs
bounded by $H$ and $H'$ are disjoint, then
$$
-2\log d_{\xi_\sharp,H_\sharp}(a,a') \leq d(H_\sharp,H) + d
 (H_\sharp,H'),
$$
with equality if and only if $d_{alg}(H,H')=0$, that is when $H$ and
$H'$ meet in one point. The corollary then follows by definition of
the inner radius.  \cqfd

\medskip
\rem Assume that $H,H'$ are circles, bounding disjoint open discs, in
the Euclidean upper halfplane, tangent to the horizontal line $\RR$ at
points $x,x'$, with Euclidean radius $r,r'>0$, and let $u=d(x,x')$ be
the Euclidean distance between $x,x'$. Then, by Pythagoras's equality
on the Euclidean right-angled triangle with base the segment between
the Euclidean centers of $H,H'$, and right-angled vertex on the
vertical line through the Euclidean center of $H$, we have
$$
rr'\leq \frac{u^2}{4}\;.
$$
Corollary \ref{coro:disjointhamenstadt} is a generalization of the
above inequality. Indeed, one may assume that $r,r'\leq \frac{1}{2}$
by homogeneity. Let $X$ be the upper halfplane model of the real
hyperbolic plane $\HH^2_\RR$, with $\xi_\sharp$ the point at infinity,
and $H_\sharp$ the horizontal line at Euclidean height $1$, so that
$\partial X=\RR\cup\{\infty\}$. The above inequality follows from
Corollary \ref{coro:disjointhamenstadt}, as $H,H'$ are horospheres in
$X$, the Hamenst\"adt distance coincides with the Euclidean distance,
and an easy computation gives $d(H_\sharp,H)=-\log(2r)$.

\bigskip
We end this section of preliminary results by the following lemma.
Let $r:[0,+\infty)\ra X$ be a geodesic ray in $X$ converging to
$\xi_\sharp$, with $r(0)$ belonging to $H_\sharp$. The {\it Busemann
function} associated to $H_\sharp$ is (see for instance
\cite{Bourdon}) the convex map $\beta_{\sharp}:X\ra \RR$ defined by
$$
\beta_{\sharp}(x)=\lim_{t\ra\infty} t-d(x,r(t))\;,
$$
so that $H_\sharp=\beta_\sharp^{-1}(\{0\})$ and $\beta_\sharp(r(t))$
tends to $+\infty$ as $t\ra +\infty$.  For any subset $A$ of $X$
define the {\it height} of $A$ with respect to $H_\sharp$ as
$$ht(A)=\sup_{x\in A} \beta_{\sharp}(x)\;.$$ 
We will say that a subset $A$ of $X$ is {\it lower} than a subset $B$
of $X$ with respect to $\xi_\sharp$ if
$$
\sup_{x\in A} \beta_\sharp(x)\leq \sup_{x\in B} \beta_\sharp(x)\;.
$$
Note that this does not depend on the chosen horosphere $H_\sharp$
centered at $\xi_\sharp$.

\begin{lemma}\label{lem:lowgeod} 
Let $H$ be a horosphere in $X$, bounding an open horoball which is
disjoint from the open horoball bounded by $H_\sharp$. Let $\eta,\eta'$
be two distinct points in $\O_{\xi_\sharp}H$. Then the geodesic
$]\eta,\eta'[$ does not intersect $H_\sharp$.
\end{lemma}

\dem 
Recall that $t\mapsto
\eta_t$ and $t\mapsto \eta'_t$ are the geodesic lines converging to
$\eta,\eta'$ respectively, with $\eta_0$ and $\eta'_0$ in $H_\sharp$.
Let $u,u'$ be points in the intersection of $H$ with
$]\xi_\sharp,\eta[,]\xi_\sharp,\eta'[$ respectively. By convexity of
the horoballs, the geodesic segment $[u,u']$ between $u$ and $u'$ is
contained in the horoball $H\!B$ bounded by $H$, hence its height is
negative. Let $t$ be big enough, so that $\eta_t$ lies between $u$ and
$\eta$ on $]\eta_\sharp,\eta[$, and similarly for $\eta'_t$. Then by
convexity of $\beta_\sharp$, the height of the geodesic segment
between $\eta_t$ and $\eta'_t$ is lower than the height of $[u,u']$,
and in particular is negative. The result follows.
\cqfd
\medskip

\section{A metric uncovering theorem}
\label{sect:abstractuncov}

In this section, we study the ``scaling'' properties in metric spaces
of families of balls, which satisfy a ``quadratic packing condition''
on the radii and the mutual distances between the centers.  The main
result, Theorem \ref{thm:abstract}, needs some elementary requirement
on the metric, as we could get a counter-example by taking discrete
spaces.  We start the section by giving a few definitions about metric
spaces.

\medskip 
If $B$ is a ball of center $x$ and radius $r$ and $\epsilon\geq 0$,
denote by $\epsilon B$ the ball of center $x$ and radius $\epsilon r$.
If $s\geq 0$ and $\B=(B_\alpha)_{\alpha\in A}$ is a family of closed
balls, define $s\B$ to be the family $(sB_\alpha)_{\alpha\in A}$ and
$s\stackrel{\circ}{\B}$ to be the family of associated open balls.

A metric space is {\it proper} if every closed ball is compact, and is
{\it geodesic} if between any two points there exists a {\it geodesic
  segment}, i.e.~a path between these two points which is an isometry
from some interval of $\RR$ into the metric space.  In any metric
space $(Y,\delta)$, such that there exists a path of finite length
between every two points, the {\it associated length distance} on $Y$
is the distance $\delta_\ell$ on $Y$ defined to be the infimum of the
lengths of paths between two points, see \cite[page 2]{GLP}. Clearly,
the length distance satisfies $\delta_\ell \geq \delta$.

A metric space {\it has extendable spheres} if for every $\epsilon>0$,
there exists $\delta$ in $]0,1[$ such that for every $r>0$ and $x,y$
in $Y$, if $(1-\delta)r\leq d(x,y)\leq r$, then there exists $y'$ in
$Y$ such that $d(x,y')=r$ and $d(y,y')\leq \epsilon \,r$. We will call
any map $\epsilon\mapsto \delta$ satisfying the above requirement a
{\it modulus of sphere extendability}.

\bexems\label{exem:preabstract}
  {\rm (i) If $M$ is a nonpositively curved complete simply
  connected Riemannian manifold, with Riemannian distance $d$, then
  $(M,d)$ has extendable spheres, with modulus $\delta(\epsilon)=
  \epsilon$. 
  
  More generally, any metric space such that through any two points
  passes at least one geodesic line (i.e.~the image of an isometry
  from $\RR$ to $Y$), has extendable spheres with modulus
  $\delta(\epsilon)= \epsilon$. This is the case for instance for the
  affine buildings with their natural metrics, and for the
  Moussong-Davis CAT$(0)$ geometric realisation of any general
  building, see for instance \cite{Dav}.
  
  But note that the Heisenberg group $\H$ (the three-dimensional simply
  connected non abelian nilpotent Lie group, see below), with its
  standard left-invariant Riemannian metric, does not have the
  property that through any two points passes at least one geodesic
  line.  Indeed, by \cite{Mar}, the only (minimizing) geodesic rays
  starting from the origin are the ones orthogonal at the origin to
  the center $[\H,\H]$ of $\H$.  Similarly, almost no (in a sense we
  do not make precise here) germ of geodesic segment for the
  Carnot-Carath\'eodory distance on $\H$ starting from the origin can be
  extended to a (minimizing) geodesic ray, see for instance \cite{BM}.
  
  \medskip (ii) The sphere $\SS^n$ for $n\geq 1$ does not have
  extendable spheres, nor does any compact metric space having at
  least two points, as the definition implies that if $(Y,d)$ has at
  least two points, then any sphere is non empty.

  \medskip (iii) Recall that the Heisenberg group $\H$ is the Lie group
  whose underlying manifold is $\CC\times\RR$, with the group law
  $$(\zeta,v)(\zeta',v')=(\zeta+\zeta', v + v' + 2\,{\rm Im}\,
  \zeta\overline{\zeta'})\;.$$
  We can endow $\H$ with the {\it Cygan distance}, which is the
  left-invariant distance $d_{\rm Cyg}$ on $\H$ such that
  $$d_{\rm Cyg}((0,0),(\zeta,v))=\sqrt[4]{|\zeta|^4 + v^2}\;,$$
  and with the {\it Carnot-Carath\'eodory distance}, which is (see
  \cite[page 161]{Gol}) the length distance $d_{\rm CC}$ associated to
  the distance $d_{\rm Cyg}$. It is also left-invariant and satisfies
  (see \cite[Theo.~4.11]{BM})
  $$d_{\rm Cyg}\leq d_{\rm CC}\leq \sqrt{\pi} \;d_{\rm Cyg}\;.$$
  
  \bprop \label{prop:Heisextspher} 
  The metric space $(\H,d_{\rm CC})$ is a proper geodesic metric space
  having extendable spheres, with modulus $\delta(\epsilon)=1-
  (1 +\epsilon^2/\pi)^{-1/2}$.  
  \eprop
  
  \dem For every $\epsilon>0$, define $\delta=1- (1 +\epsilon^2/\pi)
  ^{-1/2}$, which belongs to $]0,1[$. Take $x=(0,0)$, $r>0$ and
  $y=(\zeta,v)$ in $\H$ such that
  $$(1-\delta)r\leq d_{\rm CC}(x,y)\leq r\;.$$
  For $t>0$, let $h_t:\H\ra\H$ be the group morphism defined by
  $(\zeta,v)\mapsto(t\zeta,t^2v)$. Note that $d_{\rm
    Cyg}(h_t(u),h_t(u'))=t \,d_{\rm Cyg}(u,u')$ for every $u,u'$ in
  $\H$, hence 
  $$d_{\rm CC}(h_t(u),h_t(u'))=t \,d_{\rm CC}(u,u')\;.$$
  Let $\alpha=d_{\rm CC}(x,y)$, and define
  $y'=(\zeta',v')=(\frac{r}{\alpha}\zeta, \frac{r^2}{\alpha^2}v)$, so
  that $d_{\rm CC}(x,y')=r$. Now
  $$d_{\rm CC}(y,y')\leq \sqrt{\pi}\;d_{\rm Cyg}(y,y') =
  \sqrt{\pi}\;\big(|\zeta'-\zeta|^4 + |v'-v-2\,{\rm
    Im}\,\zeta'\overline{\zeta}\;|^2\big)^{\frac{1}{4}}
  $$
  $$= \sqrt{\pi}\;\left(\big(\frac{r}{\alpha}-1\big)^4|\zeta|^4 +
    \big(\frac{r^2}{\alpha^2}-1\big)^2 v^2\right)^{\frac{1}{4}} \leq
  \sqrt{\pi}\left(\frac{r^2}{\alpha^2}-1\right)^{\frac{1}{2}}
  \sqrt[4]{|\zeta|^4 + v^2}
  $$
  $$\leq \sqrt{\pi}\left(\frac{1}{(1-\delta)^2}-1\right)
  ^{\frac{1}{2}} d_{\rm Cyg}(x,y) \leq\epsilon \;d_{\rm CC}(x,y)\leq
  \epsilon\; r\;.
  $$
  Hence the result follows, using the homogeneity to pass from
  $x=(0,0)$ to any $x$.  \cqfd

  \medskip (iv) More generally, let $(Y,d)$ be a proper metric space.
  Assume that $(Y,d)$ {\it has dilations}, i.e.~that for every point
  $x$ in $Y$, there exists a one-parameter group
  $(h_{x,t})_{t\in\RR}$ of homeomorphisms of $Y$, such that
  $h_{x,t}\mid_{B(x,1)}$ converges uniformly to the identity map of
  $B(x,1)$ as $t$ goes to $0$, uniformly in $x$, such that
  $h_{x,t}(x)=x$, and such that $$d(h_{x,t}(y),h_{x,t}(y'))=e^t\,d(y,y')$$
  for every $t$ in $\RR$ and $y,y'$ in $Y$.  Then $(Y,d)$ has extendable
  spheres.
  
  Indeed, let $\delta>0$, $r$ in $\RR$ and $x,y$ in $Y$ such that
  $0<(1-\delta)e^r\leq d(x,y)\leq e^r$. With $\alpha=d(x,y)$, let $y'=
  h_{x,r-\log\alpha}(y)$. Then $d(x,y')=e^r$. Furthermore,
  $$d(y,y')= e^r\,d\left(h_{x,-r}(y),h_{x,r-\log \alpha}
    \big(h_{x,-r}(y)\big)\right)\;.$$
  Note that $h_{x,-r}(y)$ remains in $B(x,1)$, and that $r-\log
  \alpha$ converges to $0$ as $\delta$ goes to $0$, uniformly in $r,x$
  and $y$.  
} \eexems

\medskip 
The main result of this section is the following one.

\medskip
\begin{theorem} \label{thm:abstract}
  Let $(Y,d)$ be a proper geodesic metric space with extendable
  spheres and $0<D\leq \frac{1}{4}$. Then there exists
  $s_0=s_0(D,\delta)>0$, depending only on $D$ and on a modulus of
  sphere extendability $\delta$ of $(Y,d)$, such that the following
  holds. Let $\B= \big(B_\alpha=
  B(\xi_\alpha,r_\alpha)\big)_{\alpha\in A}$ be a finite or countable
  family of balls, having uniformly bounded radii and satisfying the
  following condition:
\begin{equation}
  0<r_\alpha r_\beta\leq  D\, d(\xi_\alpha,\xi_\beta)^2 
\label{eqn:quadratic condition}
\end{equation}
  for every distinct $\alpha,\beta$ in $A$.  
  Then, the family $s\stackrel{\circ}{\B}$ does not cover $Y$ if
  $0\leq s<s_0$. 
\end{theorem}

\begin{remark}  \label{rem:postabstract}
  (i) The quadratic packing condition \eqref{eqn:quadratic condition}
  is a generalisation of the inequality \eqref{eqn:disjhamineq} which
  is satisfied by the inner radii of the shadows of disjoint
  horoballs.
  
  \medskip (ii) Some restrictions like the quadratic packing condition
  \eqref{eqn:quadratic condition} and the assumption on an upper bound
  on the radii are necessary:
  
  Let $(Y,d)$ be the real line with the usual Euclidean metric.  Let
  $\B=\big(B(\alpha,1))_{\alpha\in \QQ}$. Clearly, the family
  $s\stackrel{\circ}{\B}$ covers the real line for all $s$. This shows
  that the boundedness of the radii is not sufficient.
  
  Let $H\!B_0$ be the horoball in $\hdr$ of Euclidean radius $1$
  tangent at $0$ to the horizontal line.  The transformation
  $Az=16z-8$ maps $H\!B_0$ to the horoball of Euclidean radius $16$
  tangent at $-8$ to the horizontal line. It is easy to check (using
  Pythagoras' theorem) that this horoball is tangent to $H\!B_0$.
  Thus, the sequence $\big(H\!B_n=A^n(H\!B_0)\big)_{n\in\ZZ}$ consists
  of horoballs with disjoint interiors. All the horoballs are tangent
  to the $A$-invariant Euclidean line that passes through the fixed
  point $8/15$ of $A$ and is tangent to $H\!B_0$. The sequence
  $\B'=(\O_\infty H\!B_n)_{n\in\ZZ}$ of the shadows of these horoballs
  is a family of intervals in $\RR$ satisfying the quadratic packing
  condition \eqref{eqn:quadratic condition}, because the horoballs are
  disjoint and by Corollary \ref{coro:disjointhamenstadt}.  However,
  it is clear that for $s>0$, the scaled collection
  $s\stackrel{\circ}{\B'}$ covers the real line.
  
  For another more arithmetic example, let $\H\!\B$ be the
  Apollonian packing of the upper halfspace: if  $H\!B_\infty$ denotes 
  the closed halfspace above the horizontal line at Euclidean height
  $1$, and $\Gamma$ is the modular group ${\rm PSL}_2(\ZZ)$,  with
  $\Ga_\infty$ the stabiliser of the point at infinity, then
  $\H\!\B=(\ga H\!B_\infty)_{\ga\in\Ga/\Ga_\infty}$. More explicitly,
  $\H\!\B=(H\!B_r)_{r\in\QQ\cup\{\infty\}}$ where $H\!B_{p/q}$,
  $p/q\in\QQ$, is the
  horoball with Euclidean center $\left(p/q ,1/(2q^2)\right)$ and
  Euclidean radius $1/(2q^2)$.   
  Let $\xi$ be a
  real number whose continued fraction expansion is unbounded. Let
  $\alpha$ be an element of ${\rm PSL}_2(\RR)$ sending $\xi$ to
  $\infty$.  Let $\B''$ be the family of shadows seen from $\infty$ of
  the horoballs in $\alpha\H\!\B$, which satisfy the quadratic packing
  condition \eqref{eqn:quadratic condition}, because the horoballs are
  disjoint. By the properties of the continued fraction expansion (see
  for instance 
  \cite[Thm.~23]{Khi}), any scaled collection
  $s\stackrel{\circ}{\B''}$ covers the real line. 
  
  \medskip (iii) If $(Y,d)$ is a proper geodesic metric space, such
  that through any two points of $Y$ passes a geodesic line, and if we
  define, for every $D$ in $]0,\frac{1}{2}[$, 
  $$
  s_0(D)=\frac{\sqrt{1+16 D}-1-4D}{4 D}\;.
  $$
  then the proof of Theorem \ref{thm:abstract} shows that the
  conclusion of Theorem \ref{thm:abstract} holds for every $s\leq
  s_0(D)$. This is a stronger result, because $D$ can be taken bigger,
  and $s_0(D)$ depends only on $D$. We have decided not to emphasize
  this point, as we do not have any example (besides the real
  hyperbolic space) where the Hamenst\"adt distance on the boundary
  minus a point of a negatively curved complete simply connected
  Riemannian manifold satisfies the property that through any two
  points passes a geodesic line.
  
  \medskip (iv) The constant $s_0$ is probably not optimal in general,
  even under the stronger assumption on $(Y,d)$ of the previous
  remark, but see Section \ref{sec:2-dim} and \ref{sec:hnr} for
  similar theorems with sharp constants in particular cases.
\end{remark}

\medskip
Before proving Theorem \ref{thm:abstract}, we start by the elementary
observation that annuli in geodesic metric spaces contain balls of
a definite size.  This is used to prove Proposition
\ref{prop:twoabstshad} which provides the induction step in the proof
of Theorem \ref{thm:abstract}.

\begin{lemma}\label{lem:canonicalball}
  Let $(Y,d)$ be a geodesic metric space. Let $\xi\in Y$ and
  $0<r_1<r_2$.  Let $p\in S(\xi,r_2)$. Then there is a
  (closed) ball $K$ of radius $(r_2-r_1)/2$, whose center is at
  distance $(r_2+r_1)/2$ from $\xi$, contained in the annulus
$$
  \overline{B(\xi,r_2)-B(\xi,r_1)}\;,
$$
  containing $p$ and meeting $B(\xi,r_1)$.
\end{lemma}

\dem 
There exists a geodesic segment in $B(\xi,r_2)$ which connects
$p$ to $\xi$.  Let $\eta$ be the point at distance $(r_2-r_1)/2$ from
$p$ on this geodesic segment.  A simple argument using the triangle
inequality shows that the ball $K=B(\eta,(r_2-r_1)/2)$ satisfies the
claim.
\cqfd
\medskip

We call any ball $K$ satisfying the conclusion of Lemma
\ref{lem:canonicalball} a {\em canonical ball} in the annulus
$\overline{B(\xi,r_2)-B(\xi,r_1)}$.

\begin{proposition}\label{prop:twoabstshad} 
  Let $(Y,d)$ be a geodesic metric space with extendable spheres, and
  $0<D\leq 1/4$. Then there exists $s_0=s_0(D,\delta)>0$, depending
  only on $D$ and on a modulus of sphere extendability $\delta$ of
  $(Y,d)$, such that the following holds. Let $\xi,\xi'$ be points in
  $Y$ and $r\ge r'>0$ such that
\begin{equation}\label{eqn:quadpacksimp}
rr'\le D\, d(\xi,\xi')^2.
\end{equation} 
If $0\leq s< s_0$, then for every canonical ball $K$ in
$\overline{B(\xi,r)- B(\xi,s r)}$, either $K$ does not meet $B(\xi',s
r')$ or there exists a canonical ball $K'$ in $\overline{B(\xi',r')-
  B(\xi',s r')}$ which is contained in $K$.
\end{proposition}

\dem Let us begin by fixing the value of $s_0$. The significance of
this particular value of $s_0$ will become apparent in the proofs of
Cases 1 and 2 below.
Consider the map $\phi:\;]0,+\infty[\;\ra \RR $ defined by
$$\phi(t)=\frac{\sqrt{4t +1} -t-1}{t}\;.$$
Note that this map is strictly decreasing on $]0,2]$ from 
$\phi(0^+)=1$ to $\phi(2)=0$, hence has values in $]0,1[$ on $]0,2[$.

Let $\delta(\cdot)$ be a modulus of sphere extendability of $(Y,d)$,
which has values in $]0,1[$.  For every $D$ in $]0,\frac{1}{4}]$, 
define
\begin{equation}\label{eqn:defsmax}
s_0=\sup_{0<\epsilon<1}\min\{\phi\big(4D(1+\epsilon)\big),
\phi\big(4D(2-\delta(\epsilon))\big)\}\;.
\end{equation}
Note that $s_0$ belongs to $]0, 1]$. Also remark that if 
$\delta(\epsilon)=\epsilon$ (as in Example \ref{exem:preabstract}
(i)), then $s_0=\phi(6D)$.

Take $s$ with $0\leq s<s_0$. Choose an $\epsilon$ in $]0,1[$ such that
$s\leq \min\{\phi\big(4D(1+\epsilon)\big),
\phi\big(4D(2-\delta(\epsilon))\big)\}$, and define $\delta=
\delta(\epsilon)$. Note that in particular, we have $s<1$.

\medskip
Let $\xi,\xi',r,r',K$ be as in the statement.  Assume that $K$
meets $B(\xi',s r')$, and let us prove the existence of $K'$ as in the
statement. As $s< 1$ the annuli
$\overline{B(\xi,r)- B(\xi,s r)}$ and $\overline{B(\xi',r')- B(\xi',s
  r')}$ are well defined.

\medskip 
Let $K={B(\eta,R)}$. By the properties of the canonical balls, we
have
$$
R=r\frac{1-s}2\;.
$$
Since $K$ meets $B(\xi',s r')$, say in a point $z$, and as $K$ is
contained in $B(\xi,r)$, we have
$$d(\xi,\xi')\leq d(\xi,z)+d(z,\xi')\leq r + sr'\leq (1+s)r\;.$$
By
Equation \eqref{eqn:quadpacksimp}, we have
$rr'\leq D\, d(\xi,\xi')^2\leq D(1+s)^2r^2$, so that
\begin{equation}
r'\leq D(1+s)^2r=\frac{2D(1+s)^2}{1-s}R\;.
\label{eqn:majorprim}
\end{equation}

\medskip\noindent{\bf Case 1 : } 
First assume that $d(\xi',\eta)\geq r'$.

\medskip Fix a geodesic segment between $\eta$ and $\xi'$. Define
$x',y'$ to be the points at distance $r',sr'$ respectively from $\xi'$
on this segment. Let $\eta'$ be the midpoint of $x',y'$ on this
segment, and $R'=r'\frac{1-s}{2}$.  Note that the points $\eta,
\eta',y'$ and $\xi'$ are in this order on this segment (see the
picture below).

\begin{figure}[ht]
\input{fig_Step2nnew.pstex_t}
\caption{}\label{fig:Step2}
\end{figure}

Define $K'=B(\eta',R')$, which is a canonical ball in
$\overline{B(\xi',r')- B(\xi',s r')}$. Let us prove that $K'$ is
contained in $K$.

Let us first show that $y'$ belongs to $K$. If this is not the case,
then $d(y',\eta)>R$.  Let $z$ be a point in $K\cap B(\xi', sr')$. The
triangle inequality gives
$$
  d(\eta,\xi')=d(\eta,y')+d(y',\xi')>R+sr'\geq
  d(\eta,z)+d(z,\xi')\geq d(\eta,\xi')\;,
$$
which is a contradiction.

Let $w$ be a point in $K'$, and let us prove that it belongs to $K$.
We have
$$
  d(w,\eta)\leq d(w,\eta')+d(\eta',\eta)\leq R'+d(\eta',\eta)=
  d(y',\eta')+d(\eta',\eta)= d(y',\eta)\leq R.
$$
Thus, $K'$ is contained in $K$.

\medskip\noindent{\bf Case 2 : } 
Assume now that $r'\geq d(\xi',\eta)\geq (1-\delta)r'$.

\medskip Then since $\delta(\cdot)$ is a modulus of sphere
extendability for $Y$, and $\delta=\delta(\epsilon)$, there exists a
point $\zeta$ in $Y$ such that $d(\zeta,\eta)\leq \epsilon \,r'$ and
$d(\zeta,\xi')=r'$. Let $K'=B(\eta',r'\frac{1-s}{2})$ be a canonical
ball in $\overline{B(\xi',r')- B(\xi',s r')}$ containing $\zeta$, so
that $\zeta$ lies on the sphere of center $\eta'$ and radius
$r'\frac{1-s}{2}$. Let us prove that $K'$ is contained in $K$.

For every $y$ in $K'$, we have
$$d(y,\eta)\leq d(y,\eta')+d(\eta',\zeta)+ d(\zeta,\eta)\leq
r'\frac{1-s}{2} + r'\frac{1-s}{2} + \epsilon
\,r'=r'(1-s+\epsilon)\;.$$
So that by Equation \eqref{eqn:majorprim}
we get
$$d(y,\eta)\leq \frac{2D(1-s+\epsilon)(1+s)^2}{1-s}R\leq
\frac{2D(1+\epsilon)(1+s)^2}{1-s}R\;.$$
Since $s\leq \phi\big(4D(1+\epsilon)\big)$, an easy
computation shows that the term on the right is at most $R$, which
proves the result.

\medskip\noindent{\bf Case 3 : } 
Finally assume that $d(\xi',\eta)\leq (1-\delta)r'$.

\medskip 
Hence for every $y$ in $B(\xi',r')$, we have, using again
Equation \eqref{eqn:majorprim},
$$d(\eta,y)\leq d(y,\xi')+d(\xi',\eta)\leq r' +(1-\delta)r' \leq
\frac{2D(2-\delta)(1+s)^2}{1-s}R\;.$$
Since $s\leq \phi\big(4D(2-\delta)\big)$, an easy computation shows
that $\frac{2D(2-\delta)(1+s)^2}{1-s}\leq 1$, so that $y$ belongs to
$K$.  Hence any canonical ball in $\overline{B(\xi',r')- B(\xi',s
  r')}$ is contained in $K$.  Note that as
$${r'} {\,}^2\leq rr'\leq D\,d(\xi,\xi')^2\leq d(\xi,\xi')^2\;,$$
we have $d(\xi,\xi')\geq r'$. Since $(Y,d)$ is a geodesic metric
space, the sphere of center $\xi'$ and radius $r'$ is non empty, hence
such a canonical ball does exists by Lemma \ref{lem:canonicalball}.
\cqfd

\medskip 
\rem
To prove Remark \ref{rem:postabstract} (iii), we make an analogous
proof as the one above, adding the assumption on $Y$ that through any
two points of $Y$ passes a geodesic line, and requiring only $D$ to be
in $]0,1/2[$. Define now $s_0=\phi(4D)$. Let $\xi,\xi',r,r',s,K$ be as
in the statement of Proposition \ref{prop:twoabstshad}. Assume that
$K$ meets $B(\xi',s r')$.  The proof of Case $1$ is unchanged. For the
two other cases, that is if $d(\xi',\eta)< r'$, define $\zeta$ to be
the point at distance $r'$ from $\xi'$ on the same side as $\eta$ on
some geodesic line $\ell$ passing through $\eta$ and $\xi'$. Let
$K'=B(\eta',r'(1-s)/2)$ be a canonical ball in $\overline{B(\xi',r')-
  B(\xi',s r')}$ containing $\zeta$.

If $\eta$ lies between $\eta'$ and $\xi'$ on $\ell$, then for every
$y$ in $K'$,
$$d(y,\eta)\leq d(y,\eta')+d(\eta',\eta)\leq r'\frac{1-s}{2}+
r'\frac{1+s}{2}=r'\;.$$

If $\eta$ lies between $\zeta$ and $\eta'$ on $\ell$, then for every
$y$ in $K'$, as both $\eta$ and $y$ are in $K'$, 
$$d(y,\eta)\leq 2r'\frac{1-s}{2}\leq r'\;.$$
Using Equation \eqref{eqn:majorprim}, we have $d(y,\eta)\leq
2D(1+s)^2R/(1-s)$.  Hence, as $s<s_0=\phi(4D)$, the points $y$
belongs to $K$, which proves the result.

Taking $D=1/4$ gives Proposition \ref{propintro} in the introduction.

\bigskip
\noindent{\bf Proof of Theorem \ref{thm:abstract}. } 
Let $Y,D$ be as in the statement. Take $s$ with $0<s < s_0$ where
$s_0=s_0(D,\delta)\leq 1$ is as in Proposition \ref{prop:twoabstshad}.
Define
$$\epsilon=1-(1+s)\sqrt{D}\;,$$
which is easily seen to be positive,
as $D\leq\frac{1}{4}$ and $s_0<1$. 

Let $\B$ be as in the statement. We may assume that $A$ is non empty.
As the balls in $\B$ have uniformly bounded radii, we can assume that
$\sup_{\alpha\in A}r_\alpha=1$, up to replacing $Y$ by
$\frac{1}{\sup_{\alpha\in A}r_\alpha}\,Y$ (which satisfies the same
assumptions as $Y$) and $\B$ by $\frac{1}{\sup_{\alpha\in
    A}r_\alpha}\,\B$, since the quadratic packing condition
\eqref{eqn:quadratic condition} is scaling invariant.

\medskip 
Let us prove the stronger statement that given $\alpha_0$ in
$A$ such that $r_{\alpha_0}\geq 1-\epsilon$, then there exists at
least one point in $B_{\alpha_0}$ that does not belong to
$\bigcup_{\alpha\in A} \,s\!\stackrel{\circ}{B}_\alpha$.

\medskip
For every $\alpha$ in $A-\{\alpha_0\}$ such that $r_\alpha>
r_{\alpha_0}$, by the quadratic packing condition, we have
$$
  d(\xi_\alpha,\xi_{\alpha_0})>\frac{1}{\sqrt{D}}(1-\epsilon)\;.
$$
If the balls $B_{\alpha_0}$ and $sB_{\alpha}$ were meeting, then
$$d(\xi_\alpha,\xi_{\alpha_0})\leq r_{\alpha_0} +sr_\alpha\leq 1
+s\;,$$
which would contradict the definition of $\epsilon$.  Hence,
up to removing such $\alpha$'s, we may assume that $r_\alpha\leq
r_{\alpha_0}=1$ for every $\alpha$ in $A$.

For every $\alpha$ in $A-\{\alpha_0\}$ such that
$d(\xi_\alpha,\xi_{\alpha_0})>3$, the balls $B_{\alpha_0}$ and
$sB_{\alpha}$ are disjoint, by the triangular inequality and as
$r_{\alpha_0},sr_\alpha$ are at most $1$. Hence, up to removing such
$\alpha$'s, we may assume that $d(\xi_\alpha,\xi_{\alpha_0})\leq 3$
for every $\alpha$ in $A$.

For every $\epsilon'>0$, by compactness of the ball $B(\xi_{\alpha_0},
3)$ and the quadratic packing condition, there are only finitely many
elements $\alpha$ in $A$ such that $r_\alpha\geq \epsilon'$. Hence we
can order the elements of $A-\{\alpha_0\}$ by
$\alpha_1,\alpha_2,\dots$ such that the radii of the balls
$B_{\alpha_i}$ are non-increasing.  Define $r_i=r_{\alpha_i}$ and
$\xi_i=\xi_{\alpha_i}$ for $i$ in $\NN$.

\medskip 
Let us construct by induction on $i$ in some initial
subsegment of $\NN$, an element $n_i$ in $\NN$ and a canonical ball
$K_i$ in $\overline{B(\xi_{n_i},r_{n_i})- B(\xi_{n_i},s r_{n_i})}$,
contained in $K_{i-1}$ for $i>0$, and which does not meet the interior
of $B(\xi_j,s r_j)$ for any $j\le n_i$.

Let $n_0=0$ and $K_0$ be any canonical ball in
$\overline{B(\xi_0,r_0)-B(\xi_0,s r_0)}$ (which exists by Lemma
\ref{lem:canonicalball}). They satisfy the requirement for $i=0$.
Assume the construction done for $i$ in $\NN$.

If $K_i$ does not meet any $B(\xi_n,s r_n)$ for $n>n_i$, then the
sequence $(n_i,K_i)$ stops. Otherwise, let $n_{i+1}$ be the first
integer $n\geq n_i+1$ for which $K_i$ meets $B(\xi_n,sr_n)$.
Proposition \ref{prop:twoabstshad} implies that there exists a
canonical ball $K_{i+1}$ in
$\overline{B(\xi_{n_{i+1}},r_{n_{i+1}})-B(\xi_{n_{i+1}},s
  r_{n_{i+1}})}$ which is contained in $K_{i}$. The couple
$(n_{i+1},K_{i+1})$ satisfies the requirement.

Clearly, any point in the nonempty intersection $\bigcap_i K_i$ of the
nested nonempty compact subsets $K_i$ is in the complement of
$\bigcup_{\alpha\in A} s \stackrel{\circ}{B}_\alpha$.
\cqfd
\medskip

\begin{remark}\label{rem:two_balls} 
  At the beginning of the induction in the proof of Theorem
  \ref{thm:abstract}, we fixed a canonical ball $K_0$. We obtained a
  point in the complement of $s\!\stackrel{\circ}{\B}$ which is
  contained in $K_0$. Assume that in $(Y,d)$ the {\it spheres have
    antipodal points}, i.e.~for every $y$ in $Y$ and $r>0$, there
  exists $p,q$ in the sphere $S(y,r)$ such that $d(p,q)=2r$. This is
  the case if through any point in $Y$ passes a geodesic line, in
  particular if $Y$ is a Carnot group endowed with its
  Carnot-Carath\'eodory metric (see Remark \ref{rem:CCantipod}). Then
  take antipodal points $p,p'$ on the sphere of center
  $\xi_{\alpha_0}$ and radius $r_{\alpha_0}$.  If $K_0,K'_0$ are
  canonical balls in $\overline{ B(\xi_{\alpha_0},r_{\alpha_0})-
    B(\xi_{\alpha_0},sr_{\alpha_0})}$ containing $p,p'$ respectively
  (which exist by Lemma \ref{lem:canonicalball}), then $K_0$ and
  $K'_0$ are disjoint (they are even at distance at least
  $sr_{\alpha_0}>0$, by the triangular inequality). Thus, in the
  conclusion of Theorem \ref{thm:abstract}, we can obtain that if
  $s<s_0$, then there exist not only one, but at least two points
  outside $s\!\stackrel{\circ}{\B}$. We will use this fact in the next
  Section (see Remark \ref{rem:doublerlamise}).
\end{remark}

\section{Unclouding the sky of negatively curved manifolds} 
\label{sec:spaboundsha}

In this section, we consider a proper 
CAT$(-1)$ 
space
$X$, a point $\xi_\sharp$ in $\partial X$, and a horoball $H_\sharp$
centered at $\xi_\sharp$. We denote by
$d_\infty=d_{\xi_\sharp,H_\sharp}$ the Hamenst\"adt distance on
$\partial X-\{\xi_\sharp\}$ and by $d_\ell$ its associated length
distance.  We assume that the following two conditions on
$(X,\xi_\sharp)$ hold:
\begin{itemize}
\item[(i)]
the metrics $d_\infty$ and $d_\ell$ are equivalent, i.e.~there
exists $C\geq 1$ such that
\begin{equation}
  d_\infty\le d_\ell\le Cd_\infty. \label{eqn:equivalence}
\end{equation}
\item[(ii)] the metric space $(\partial X-\{\xi_\sharp\},d_\ell)$ has
extendable spheres, with modulus $\delta(\cdot)$.
\end{itemize}
The inequalities \eqref{eqn:equivalence} in particular imply that
$d_\ell$ is finite, so that $(\partial X-\{\xi_\sharp\},d_\ell)$ is a
proper geodesic metric space. Note that these two properties do not
depend on $H_\sharp$. By Remark \ref{rem:hamlocequiv}, if 
$H'_\sharp$ is another horosphere in $X$, centered at the point at infinity
$\xi'_\sharp$, then the  Hamenst\"adt distance
$d'_\infty=d_{\xi'_\sharp,H'_\sharp}$ and its associated length distance
$d'_\ell$ also satisfy the first property (i) locally, and locally
$$\lambda e^{-8} d_\ell\leq d'_\ell\leq\lambda e^8 d_\ell\;$$
for some
$\lambda>0$, so that $(\partial X-\{\xi'_\sharp\},d'_\ell)$ also has
extendable spheres restricted (in a sense we won't make precise) to a
compact subset of $\partial X-\{\xi_\sharp,\xi'_\sharp\}$.

\medskip The class of spaces for which (i) and (ii) hold includes
(besides the real hyperbolic 
$n$-space $\HH_\RR^n$, with the constant $C=1$, as we have already
seen) the complex hyperbolic $n$-space $\HH_\CC^n$, with the
constant $C=\sqrt{\pi}$. To prove this, we use the upper halfspace
model of $\HH_\CC^n$ with $\xi_\sharp$ the point at infinity, where
the space $\partial \HH_\CC^n-\{\infty\}$ is identified with the
Heisenberg group $\H$. Then the Hamenst\"adt distance on $\partial
\HH_\CC^n-\{\infty\}$ is (up to a constant factor) the Cygan metric
$d_{\rm Cyg}$ on $\H$ (see \cite[page 219]{HP4}). We have already said
that the Carnot-Carath\'eodory metric $d_{CC}$ on $\H$ is the length
distance associated to the Cygan distance on $\H$, and that the Cygan
and Carnot-Carath\'eodory distances satisfy $d_{\rm Cyg}\leq d_{CC}\leq
\sqrt{\pi}d_{\rm Cyg}$. We proved in Proposition
\ref{prop:Heisextspher} that the Carnot-Carath\'eodory distance on $\H$
has extendable spheres with modulus $\delta(\epsilon)=1-
  (1 +\epsilon^2/\pi)^{-1/2}$. Hence the two properties (i), (ii)
above are satisfied for $(\HH_\CC^n,\infty)$ (and any point at
infinity $\xi_\sharp$ instead of $\infty$ by homogeneity).

More generally, the two conditions (i) and (ii) are true if $X$ is any
negatively curved symmetric space, see \cite[Lemmas 2.1, 2.2]{Ham}
(where the distance defined in this paper is not exactly the
Hamenst\"adt distance, but is equivalent to it, see \cite{HP3}), that
is, besides the above cases, for the quaternionic hyperbolic spaces and
the octonionic hyperbolic plane.

\medskip
There are infinitely many more examples. Let us first recall the basic
structural properties of a connected homogeneous negatively curved
Riemannian manifold $M$, which are essentially due to \cite{Kob,Wolf,Hei} (see
also \cite{Ale,AW1,AW2}).

Such a manifold $M$ is simply connected and there exists a simply
transitive solvable Lie group $S$ of isometries of $M$
\cite[Prop.~1]{Hei}. We identify $M$ and $S$ from now on, so that $S$
carries a left-invariant negatively curved Riemannian metric.  Let
$N=[S,S]$, and let ${\mathfrak S}, {\mathfrak N}$ be the Lie algebras of
$S$ and $N$. If $X$ is a unit length vector in ${\mathfrak S}$ which is
orthogonal to ${\mathfrak N}$, then ${\mathfrak S}=\RR X \oplus
{\mathfrak N}$ and (up to changing $X$ to $-X$) the eigenvalues of the
derivation $A= {\rm ad~} X\mid_{\mathfrak N}$ of ${\mathfrak N}$ have
positive real parts (see \cite[Prop.~2]{Hei} and
\cite[Sect.~1.3]{EH}). Note that $S$ and $A$ are not uniquely
determined by $M$, but ${\mathfrak N}$ is (up to isomorphism) (see
\cite{AW1,AW2} who give the precise relation between two such $S$'s,
and see also \cite[Sect.~1.8]{EH}). Conversely, if ${\mathfrak N}$ is
a nilpotent Lie algebra endowed with a derivation $A$ whose
eigenvalues have positive real parts, then the simply connected Lie
group, whose Lie algebra is the one-dimensional extension of
${\mathfrak N}$ constructed using $A$, carries a left-invariant
negatively curved Riemannian metric (see \cite[Theo.~3]{Hei}). This
metric is not unique in general (see for instance \cite{EH} for
precise examples).

\medskip
Now, let ${\mathfrak N}_\lambda$ be the 
eigenspace of $A$
associated to a complex number $\lambda$ (which is trivial if
$\lambda$ is not an eigenvalue). 
Note that $[{\mathfrak
N}_\lambda,{\mathfrak N}_\mu]$ is contained in ${\mathfrak
N}_{\lambda+\mu}$, as $A$ is a derivation.

We say that $M$ is a homogeneous negatively curved Riemannian manifold
{\it of Carnot type} if
${\mathfrak N}_1$ generates ${\mathfrak N}$ as a Lie algebra.  In
particular, $N$ is a {\it Carnot group} (see for instance
\cite{FS,Pan2,Gro3}), which means that ${\mathfrak N}$ is graduated by
${\mathfrak N}=\bigoplus_{i=1}^k{\mathfrak N}_i$ (for some $k$) as a
Lie algebra, with ${\mathfrak N}_1$ generating ${\mathfrak N}$ as a
Lie algebra.  We will make some comments on this definition after the
proof of the following proposition.

\bprop \label{prop:homonegcurvriemancarnok}
Any homogeneous negatively curved Riemannian manifold of Carnot
type satisfies the above two properties (i) and (ii), for some point
at infinity.  
\eprop

\dem Let $\exp :{\mathfrak S}\ra S$ be the exponential map.  Note that
by \cite{Hei} (see also \cite{EH}), the curve $t\mapsto \exp tX$ is a
(unit speed) geodesic line passing through the identity $e$ of $S$ at
time $t=0$. Furthermore, $N$ is the horosphere containing $e$ centered
at the point at infinity $\xi_\sharp$ of $t\mapsto \exp tX$. Identify
the punctured boundary at infinity $\partial S-\{\xi_\sharp\}$ with
$N$ by the map sending the point $g$ in $N$ to the endpoint of the
geodesic line starting from $\xi_\sharp$ and passing through $g$ at
time $t=0$.  Note that this geodesic line is $t\mapsto g\exp(
-tX)$. Hence the isometry of $S$ which is the left translation by
$\exp(tX)$, fixes the point $\xi_\sharp$ and extends to the punctured
boundary at infinity, identified to $N$ as above, by the map
$\psi^t:g\mapsto \exp(tX)g\exp( -tX)$.

As $\exp(tX) N$ is the horosphere centered at $\xi_\sharp$ containing
$\exp tX$, its algebraic distance to $N$ is exactly $t$. By the
formulas \eqref{eqn:invadistham} and \eqref{eqn:disthamvarihoro} in
Section \ref{sec:geomhoro}, for every $a,b$ in $\partial
S-\{\xi_\sharp\}$, we have
$$d_{\xi_\sharp,N}(\psi^t(a),\psi^t(b)) =
d_{\xi_\sharp,\exp(-tX)N}(a,b)=e^t \;d_{\xi_\sharp,N}(a,b)\;.$$
Hence $(\psi^t)_{t\in\RR}$ is a one-parameter group of dilations for
the Hamenst\"adt distance $d_{\xi_\infty,N}$.

\medskip 
Recall that the Carnot-Carath\'eodory distance $d_{CC}$ on $N$
is defined as follows (see for instance \cite{Gro3}).  Endow
${\mathfrak N}_1$ with a Euclidean norm 
(for instance the one induced by the Riemannian metric on $S$).
Consider the left-invariant Euclidean vector sub-bundle $\Delta$ of
$TN$ defined by the Euclidean subspace ${\mathfrak N}_1$ of $T_eN$.
Then $d_{CC}(x,y)$ is the lower bound of the Euclidean lengths of the
smooth paths between $x$ and $y$ in $N$ that are tangent to $\Delta$
at every point.  

\medskip We make an independent remark, which will be used in Section
\ref{sec:bi-infinite}. 

\brema \label{rem:CCantipod} In a Carnot group endowed with its
Carnot-Carath\'eodory distance, the spheres have antipodal points (in
the sense of Remark \ref{rem:two_balls}). 
\erema

\dem
Indeed, the map $t\mapsto exp (tX)$ for $X$ in ${\mathfrak N}_1$
is a geodesic line (parametrized proportionally to arc length) through
$e$, and then one uses the homogeneity of $N$.  
\cqfd

\medskip Recall that
$$\exp tX \exp Y \exp -tX= \exp \left( e^{t\;{\rm ad}\;X}(Y)\right)\;,$$
and that $\exp:{\mathfrak N}\ra N$ is a diffeomorphism.  As $\psi^t$ is
a group morphism of $N$ and as $e^{t{\rm ~ad~}X}$ is the homothety of
ratio $e^t$ on ${\mathfrak N}_1$, the map $\psi^t$ acts on the Euclidean
sub-bundle $\Delta$ by a dilation of ratio $e^t$ from $\Delta_g$ to
$\Delta_{\psi^t(g)}$. Hence $(\psi^t)_{t\in\RR}$ is also a
one-parameter group of dilations for the Carnot-Carath\'eodory distance
$d_{CC}$.

Since the Carnot-Carath\'eodory and Hamenst\"adt distances induce the
same topology, they are equivalent on the spheres of center $e$ and
radius $1$ in $N$. As $(\psi^t)_{t\in\RR}$ is a one-parameter group of
dilations for both of them, they are also equivalent on every sphere
of center $e$. As they are both invariant by left-translation under
$N$, these two distances are equivalent.

As $d_{\xi_\sharp,N}$ is equivalent to $d_{CC}$ which is a length
distance, it follows easily that the length distance $d_\ell$
associated to $d_{\xi_\sharp,N}$ is finite, and equivalent to
$d_{CC}$, hence to $d_{\xi_\sharp,N}$.

By Example \ref{exem:preabstract} (iv), since $(\partial
S-\{\xi_\sharp\},d_{\xi_\sharp,N})$ has dilations, it has extendable
spheres. This proves the result.  
\cqfd

\medskip
\begin{remark}  \label{rem:CCnece}
  (i) Many authors have considered homogeneous negatively curved
  Riemannian manifolds $M$ of Carnot type (see \cite{Gro3}). For
  instance, it is proved in \cite{HK} that the local and global
  definitions of quasi-conformal maps on the boundary of $M$ coincide.
  The negatively curved $3$-steps Carnot solvmanifolds of \cite{EH}
  belong to this class (they are precisely the ones such that
  ${\mathfrak N}={\mathfrak N}_1\oplus{\mathfrak N}_2$ and ${\mathfrak
  N}_2$ is the full center of ${\mathfrak N}$). 
  
\medskip (ii) If $M$ is a negatively curved symmetric space different
from the real hyperbolic space, then $M$ is a negatively curved
$3$-step Carnot solvmanifold, such that furthermore the center of
${\mathfrak N}$ has dimension $1,3$ or $7$. There exist infinitely
many pairwise non isomorphic Carnot groups with order of nilpotency at
least $3$ (if ${\mathfrak G}$ is any finite dimensional nilpotent Lie
algebra, with ${\mathfrak G_0}={\mathfrak G}$ and ${\mathfrak
  G_{i+1}}=[{\mathfrak G},{\mathfrak G_i}]$, then the graded Lie
algebra associated to ${\mathfrak G}$, which is
$\bigoplus_{i}{\mathfrak G_{i}}/{\mathfrak G_{i+1}}$ with the obvious
Lie bracket, is the Lie algebra of a Carnot group). There exist
infinitely many Carnot groups with order of nilpotency $2$ whose
center has dimension at least $8$ (see for instance \cite{Kap}). Each
Lie algebra ${\mathfrak N}$ of a Carnot group carries a derivation $A$
whose eigenvalues have positive real parts, defined, if ${\mathfrak
  N}=\oplus_{i=1}^k{\mathfrak N}_i$ with $[{\mathfrak N}_i,{\mathfrak
  N}_j]\subset {\mathfrak N}_{i+j}$, by $A\mid_{{\mathfrak
    N}_j}=j\;{\rm id}_{{\mathfrak N}_j}$. As seen above, we get
infinitely many (pairwise non homothetic) non symmetric homogeneous
negatively curved Riemannian manifolds.

  \medskip (iii) Let $M$ be a homogeneous negatively curved
  Riemannian manifold.  If we assume that the curvature of $M$ is
  normalized so that its maximum is $-1$, and that $A$ is
   diagnalisable over $\RR$, then the assumption (in the definition of
  being of Carnot type) that ${\mathfrak N}_1$ generates ${\mathfrak
  N}$ as a Lie algebra is in fact necessary.

To prove this, we start with the following lemma, whose proof is
contained in the proof of Lemma 1, page 486 of \cite{EH}.

\blemm (Eberlein-Heber) 
If the curvature of $M$ is at at most $-1$,
then the real part of any eigenvalue $\lambda$ of $A$ is at least $1$.
\elemm

\dem As $A$ is a linear map of the Euclidean space ${\mathfrak N}$, we
can consider the symmetric and antisymmetric part
$D_0=\frac{1}{2}(A+A^*)$, $S_0=\frac{1}{2}(A-A^*)$ of $A$. Recall that
$D_0$ is symmetric and  positive definite \cite[Prop2]{Hei}.  Define
$N_0=D_0^2 + D_0S_0-S_0D_0$. Recall that (see for instance \cite{Hei})
if $Y$ is a unit vector in ${\mathfrak N}$, then the curvature in $S$
of the tangent $2$-plane generated by $X$ (where $X$ is as above) and
$Y$ is
$$K(X,Y)=-\langle N_0Y,Y\rangle\;.$$ As $N_0$ is symmetric, hence
diagonalisable, and since $K\leq -1$, the eigenvalues of $N_0$ are at
least $1$.  Let $(Y_i)_{i\in I}$ be an orthonormal basis of
eigenvectors of $D_0$, with (positive) associated eigenvalues
$(\lambda_i)_{i\in I}$.  Note that
$$\lambda_i^2=\langle D_0^2Y_i,Y_i\rangle=\langle
N_0Y_i,Y_i\rangle=-K(X,Y_i)\geq 1\;.$$
Hence $\lambda_i\geq 1$ for every $i$. In particular, $\langle
D_0\zeta,\zeta\rangle\geq ||\zeta||^2$ for every $\zeta$ in
${\mathfrak N}$.  Now let $\xi=\xi'+i\xi''$ be a complex eigenvector
for the complex eigenvalue $\lambda=\lambda'+i\lambda''$ of $A$. We
have, by an easy computation,
$$\lambda'(||\xi'||^2+ ||\xi''||^2)=\langle A\xi',\xi'\rangle+\langle
A\xi'',\xi''\rangle= \langle D_0\xi',\xi'\rangle+\langle
D_0\xi'',\xi''\rangle\geq ||\xi'||^2+ ||\xi''||^2\;.$$
Therefore $\lambda'\geq 1$.  
\cqfd

\medskip 
Now assume that the curvature of $M$ has maximum $-1$, and that $A$ is
diagonalisable over $\RR$. Identify $N$ with ${\mathfrak N}$ by the
exponential map.  Let $\Lambda$ be the set of eigenvalues of $A$.
Consider the map $d'$ from $N\times N$ to $[0,+\infty)$, which is
left-invariant by the diagonal action of $N$ on $N\times N$, such that
for $Y=\sum_{\lambda\in\Lambda} Y_\lambda$ in ${\mathfrak
N}=\bigoplus_{\lambda\in\Lambda}{\mathfrak N}_\lambda$, we have
$$d'(0,Y)=\sum_{\lambda\in\Lambda}
||Y_\lambda||^\frac{1}{\lambda}\;.$$ Note that the map $d'$ (which is
a priori not a distance) is continuous and does not vanish outside
the diagonal of $N\times N$, and that the map $\psi^t$ is a dilation
of ratio $e^t$ for the map $d'$.  Hence there is a constant $c>0$ such
that, for every $x,y$ in $N$,
$$\frac{1}{c}d'(x,y)\leq d_{\xi_\sharp,N}(x,y)\leq
cd'(x,y)\;.$$ Therefore, to check if a path has finite length or
not, we may replace the Hamenst\"adt distance by $d'$.  As the
eigenvalues of $A$ are at least $1$ and $||\;||$ is an Euclidean norm
on $N$, the only smooth paths that have finite length are the ones
that are tangent at every point to the left-invariant vector
sub-bundle $\Delta$ of $TN$ defined by the subspace ${\mathfrak N}_1$
of $T_eN$. In order for every pair of points to be joined by a finite
length curve, it is hence necessary that ${\mathfrak N}_1$ generates
${\mathfrak N}$ as a Lie algebra. 
\end{remark}

For every horosphere $H$ in a CAT$(-1)$ metric space $X$, and every
$t\geq 0$, define $H(t)$ to be the horosphere at distance $t$ from
$H$, contained in the horoball bounded by $H$.

\begin{theorem}\label{thm:bounded2}  
  Let $X$ be a proper CAT$(-1)$ space, $\xi_\sharp$ a point in
  $\partial X$, and $H_\sharp$ a horoball centered at $\xi_\sharp$.
  Assume that $(X,\xi_\sharp)$ satisfies the conditions (i) and (ii)
  as in the beginning of Section \ref{sec:spaboundsha}.
  There exists $t_0>0$ such that for every sequence $(H_n)_{n\in\NN}$
  of horospheres bounding open horoballs that are pairwise disjoint
  and disjoint from the open horoball bounded by $H_\sharp$, if $t>
  t_0$, then there exists at least one geodesic line starting at
  $\xi_\sharp$ which avoids $H_n(t)$ for every $n$ in $\NN$.
\end{theorem}

\dem We use the notations of the beginning of the section. Define
$$t_0=t_0(C,\delta)= -\log\frac{s_0(\frac{1}{4},\delta)}{2c_{\rm
    max}\,C}\;,
$$
where $s_0(\epsilon,\delta)$ is as in Theorem \ref{thm:abstract},
and $c_{\rm max}$ as in Remark \ref{rem:ombreboulegene} (and $c_{\rm
  max}=2^{-\frac{1}{a}}$ if $X$ is a Riemannian manifold with
sectional curvature $-a^2\leq K\leq -1$, with $a\geq 1$).

Let $\xi_n$ be the point at infinity of $H_n$ and
$$r_n=r_i(H_n)=\frac{1}{2}e^{-d(H_\sharp,H_n)}>0$$ be the inner radius
of the shadow of $H_n$. Consider the family of balls
$(B_n=B_{d_\ell}(\xi_n,r_n))_{n\in\NN}$ in $(\partial
X-\{\xi_\sharp\},d_\ell)$. By Corollary \ref{coro:disjointhamenstadt},
with $D=\frac{1}{4}$, we have, if $n\neq m$,
$$0<r_nr_m\leq D d_\infty(\xi_n,\xi_m)^2\leq D
d_\ell(\xi_n,\xi_m)^2\;. $$
Note that the radii $(r_n)_{n\in\NN}$ are
uniformly bounded (by $\frac{1}{2}$), and recall that the proper
geodesic metric space $(\partial X-\{\xi_\sharp\},d_\ell)$ has
extendable spheres with modulus $\delta$.  Hence, by Theorem
\ref{thm:abstract}, for every $s<s_0(\frac{1}{4},\delta)$, there
exists a point $\xi_s$ in $\partial X-\{\xi_\sharp\}$ which does not
belong to $\bigcup_{n\in\NN} sB_n$.  By Remark
\ref{rem:ombreboulegene} (and Corollary \ref{coro:ballshadsize} if $X$
is a Riemannian manifold with sectional curvature $-a^2\leq K\leq -1$,
with $a\geq 1$), as $d(H_\sharp,H_n(t))=t+ d(H_\sharp,H_n)$, we have
$$
\O_{\xi_\sharp}H_n(t) \subset
B_{d_{\infty}}(\xi_n,2c_{\rm max}\,e^{-t}r_n)\subset
B_{d_{\ell}}(\xi_n,2c_{\rm max}\,e^{-t}Cr_n)\;.
$$
Let $s(t)=2c_{\rm max}\,e^{-t}C$, so that if $t>t_0$, then
$s(t)<s_0(\frac{1}{4},\delta)$. Therefore the point $\xi_{s(t)}$ in
$\partial X-\{\xi_\sharp\}$ does not belong to $\bigcup_{n\in\NN}
\O_{\xi_\sharp}H_n(t)$, which proves the result.  
\cqfd

\bigskip
Note that if $X$ is a complete simply connected
Riemannian manifold with sectional curvature $-a^2\leq K\leq -1$ (with
$a\geq 1$), if the Hamenst\"adt distance $d_\infty$ is geodesic and if
through any two points of $\partial X-\{\xi_\sharp\}$ passes a
geodesic line for $d_\infty$, then, by Remark
\ref{rem:postabstract} (iii), the minimal value $t_0$ is
$$t_0= -\log\frac{\sqrt 5-2}{2^{1-\frac 1a}}\;.
$$
In the case $X$ is the real hyperbolic $n$-space $\hnr$, we get
$t_0=-\log(\sqrt 5-2)\approx 1.44$. In Section \ref{sec:hnr}, we will
show how to improve this, by using the fact that the boundary of
$\hnr$ minus a point is the $(n-1)$-dimensional Euclidean space to
prove a stronger version of Theorem \ref{thm:abstract}.
If $X=\hdc$, a straightforward calculation using Proposition
\ref{prop:Heisextspher} yields $t_0\approx 4.9157$.

\begin{remark}\label{rem:doublerlamise}
  Let $(X,\xi_\sharp,H_\sharp)$ be as in Theorem \ref{thm:bounded2}.
  Assume that furthermore the spheres in $(\partial
  X-\{\xi_\sharp\},d_\ell)$ have antipodal points (in the sense of
  Remark \ref{rem:two_balls}). Then it follows from Remark
  \ref{rem:two_balls} that there exists $t_0>0$ such that for every
  sequence $(H_n)_{n\in\NN}$ of horospheres bounding open horoballs
  that are pairwise disjoint and disjoint from the open horoball
  bounded by $H_\sharp$, if $t> t_0$, then there exist at least two
  geodesic lines starting at $\xi_\sharp$ which avoid $H_n(t)$ for
  every $n$ in $\NN$. We will use this fact in Section
  \ref{sec:bi-infinite}.
\end{remark}

\begin{remark}\label{rem:finquatre}
  (1) Let $T$ be a locally finite metric tree, with degrees at least
  $3$, and denote by $\ell_{\rm max}$ the upper bound on the length of
  the edges.  Recall that $T$ is a proper CAT$(-1)$ space (see for
  instance \cite{BH}). We claim that for every sequence
  $(H_n)_{n\in\NN}$ of horospheres in $T$, bounding open horoballs
  that are pairwise disjoint, for every point $x_{0}$ not contained in
  an open horoball bounded by some $H_n$, if $t> \ell_{\rm max}$, then
  there exist at least two geodesic rays in $T$ starting at $x_0$
  which avoid $H_n(t)$ for every $n$ in $\NN$.
  
  \medskip\dem Fix a geodesic ray starting from $x_0$. Follow this ray
  until it reaches its first intersection point $y_1$ with $\bigcup _n
  H_n$ or otherwise till infinity. In the first case, let $H_{i_1}$ be
  a horosphere containing $y_1$. Follow a segment starting at $y_1$
  and entering in $H_{i_1}$, and when arriving at the first vertex
  $v_1$ of $T$, follow a segment starting at $v_1$ that does not point
  towards the point at infinity of $H_{i_1}$ or back towards
  $x_0$. Let $x_1$ be the 
  intersection point with $H_{i_1}$ after $v_1$. Fix a geodesic ray
  starting at $x_1$ that does not enter $H_{i_1}$. Iterating this
  construction, we get a geodesic ray starting from $x_0$ which avoids
  $H_n(t)$ for every $n$ in $\NN$ and $t> \ell_{\rm max}$. This gives
  our first geodesic ray.  Using the fact that the degree at $y_1$ is
  at least $3$, we easily get a second geodesic ray. \cqfd

\medskip
(2) The above result is sharp, as shown by the $3$-regular tree $T_3$
whose edge lengths are equal to some constant $N$, with a family
$(H_n)_{n\in\NN}$ of horospheres bounding disjoint open horoballs,
such that the union of the (closed) horoballs bounded by the $H_n$'s
covers $T_3$. The existence of such a family is easy, but see also 
\cite{BL,Pau} for algebraic constructions.

\medskip 
(3) The above sharpness result (where $N$ can be taken
arbitrarily large) shows that it is not possible to generalize Theorem
\ref{thm:bounded2} to the full collection of CAT$(-1)$ metric spaces,
even if we assume that $\partial X$ has no isolated point to avoid
trivial examples (as for $X=\RR$).  To be more precise, there exists
no constant $c>0$ such that for every proper geodesic CAT$(-1)$ metric
space $X$, for every horosphere $H_\sharp$ centered at $\xi_\sharp$,
for every every sequence $(H_n)_{n\in\NN}$ of horospheres in $X$,
bounding open horoballs that are pairwise disjoint and disjoint from
the open horoball bounded by $H_\sharp$, if $t> c$, then there exists
at least one geodesic ray starting at $\xi_\sharp$ which avoids $H_n(t)$ for
every $n$ in $\NN$. Some conditions on $X$, as for instance (i) and
(ii), (which in particular imply that $\partial X$ is path-connected,
which excludes the case of trees for instance) are hence necessary.
\end{remark}

\section{Unclouding the sky of two-dimensional spaces}
\label{sec:2-dim}

In this section, we prove a stronger version of Theorem
\ref{thm:bounded2} in the two-dimensional manifold case. Assume in the
whole section that $M$ is any two-dimensional complete simply
connected Riemannian manifold with 
curvature $-a^2\leq K\leq -1$ (for some $a\geq 1$). In this situation,
the boundary of $M$ is homeomorphic to a circle, which allows us to
use techniques not available in higher-dimensional situations. Note
that in most two-dimensional cases there are no curves of finite
length between two distinct points in the boundary. Thus, Theorem
\ref{thm:bounded2} does not apply in general.

Define $t_1=t_1(a)$ by
$$
e^{-t_{1}}=\left\{\begin{array}{ll}
  2^{2/a}\left(\sqrt{1+2^{1-1/a}}-1-2^{-1-1/a}\right)&
  {\rm if~} a\leq 2\\
  \min\left\{ 1-2^{-2/a}\;,\;
  2^{2/a}\left(\sqrt{1+2^{1-1/a}}-1-2^{-1-1/a}\right)\right\}&
{\rm if~} a\geq 2\;.
\end{array}\right.
$$
Note that this is well defined, as the right terms are positive, and
that $t_1(a)$ converges to $+\infty$ as $a$ converges to $+\infty$.

\begin{theorem}\label{thm:2dimbounded} 
  With $M$ as above, let $\xi_\sharp$ be a point in $\partial M$, and
  $H_\sharp$ be a horosphere in $M$ centered at $\xi_\sharp$.  For
  every sequence $(H_n)_{n\in\NN}$ of horospheres bounding open
  horoballs that are pairwise disjoint and disjoint from the open
  horoball bounded by $H_\sharp$, if $t> t_1$, then there exist at
  least two geodesic lines starting at $\xi_\sharp$ which avoid
  $H_n(t)$ for every $n$ in $\NN$.
\end{theorem}

Before giving the proof, we first fix some notations. The shadows will
be taken with respect to the point at infinity $\xi_\sharp$. We denote
by $d_\infty$ the Hamenst\"adt distance $d_{\xi_\sharp,H_\sharp}$, and
the balls will be taken with respect to this distance.  For
$\eta,\eta'$ in $\partial M-\{\xi_\sharp\}$, we denote by
$[\eta,\eta']_\infty$ the closure of the connected component of
$\partial M-\{\eta,\eta'\}$ which does not contain $\xi_\sharp$, and
similarly for half-open intervals.

If $H$ is a horosphere centered at a point at infinity $\xi$
different from $\xi_\sharp$, then by convexity, there exist two and
only two geodesics lines starting from $\xi_\sharp$ and tangent to $H$
(at time $t=0$). If $\xi_+,\xi_-$ are the points at infinity of these
two geodesic lines, then by connectedness, we have $\O
\,H=[\xi_+,\xi_-]_\infty$. In particular, for any $t>0$, the closure
$\overline{\O\,H-\O\,H(t)}$ is the union of exactly two (topological)
compact intervals. Note that $\overline{B(\xi,r)-B(\xi,r')}$ for $r'<r$
  may a priori have more than two components. We claim that
\begin{equation}\label{eqn:bonneborne}
r_i(H)\leq d_\infty(\xi,\xi_\pm)\leq r_e(H)\;.
\end{equation}
The upper bound follows from the right inclusion in the statement of
Corollary \ref{coro:ballshadsize}.  The lower bound follows from its
proof (take $B'=B$ in the paragraph around Figure
\ref{fig:inclusionombreboule}, compare with $\HH_{-1}^2$ instead of
$\HH_{-a^2}^2$, which reverses the comparison inequality, and take
limits to pass from balls to horoballs).

\medskip
The main step of the proof of Theorem \ref{thm:2dimbounded} is the
following proposition.

\begin{proposition}\label{prop:2dimhoroballs} 
  Let $M$, $\xi_\sharp$ and $H_\sharp$ be as in Theorem
  \ref{thm:2dimbounded}.  Let $H$ and $H'$ be horospheres bounding
  open horoballs that are disjoint and disjoint from the open horoball
  bounded by $H_\sharp$, with $r_i(H')\leq r_i(H)$.  Let $t\geq
  t_{1}(a)$.  For every connected component $K$ of
  $\overline{\O\,H-\O\,H(t)}$, the following dichotomy holds~:
\begin{itemize}
\item either $K$ does not meet $\O\,H'(t)$
\item or there exists a connected component $K'$ in
  $\overline{\O\,H'-\O\,H'(t)}$, which is contained in $K$.
\end{itemize}
\end{proposition}

\dem Let $H,H',t,K$ be as in the statement. Up to replacing $H_\sharp$
by the horosphere centered at $\xi_\sharp$ and tangent to $H$, we may
assume that $H$ and $H_\sharp$ are tangent, which is equivalent to
assuming that $r_i(H)=1/2$.  Let $\xi,\xi'$ be the points at infinity
of $H,H'$ respectively.  Orient the topological line $\partial
M-\{\xi_\sharp\}$ so that $K$ is on the right side of $\xi$, and let
$K=[k_1,k_2]_\infty$ with $\xi,k_1,k_2$ in this order on $\partial
M-\{\xi_\sharp\}$.

\begin{figure}[ht]
\input{fig_dimde.pstex_t}
\caption{}\label{fig:dimde}
\end{figure}

{\bf Case 1:} 
Assume that $\xi'$ belongs to $K$.  We will show the stronger
conclusion that at least one of the two components of $\O\,H'-\{\xi'\}$
is contained in $K$. This follows if
$$
r_e(H')<\max\{d_{\infty}(\xi',k_1), d_{\infty}(\xi',k_2)\}.   
$$
Indeed, recall that $\O\,H'$ is contained in $B(\xi',r_e(H'))$ by
Corollary \ref{coro:ballshadsize}. Hence, if for $i=1$ (resp.~$i=2$)
we have $r_e(H')< d_{\infty}(\xi',k_i)$, then the left (resp.~right)
connected component of $\O\,H'-\{\xi'\}$ is contained in
$[\xi',k_i[_\infty$ hence in $K$. Note that in constant curvature, and
more generally if the Hamenst\"adt distance to a given point is
monotonous on the line $\partial M-\{\xi_\sharp\}$, then we may
replace the strict inequality in the displayed formula above by a
large one.

Define $u=d_{\infty}(\xi,\xi')$.  Corollary
\ref{coro:disjointhamenstadt} implies that
$$
r_e(H')=2^{1-\frac{1}{a}}r_i(H')=2^{-\frac{1}{a}}4r_i(H')r_i(H)\leq
2^{-\frac 1a}u^2.
$$
On the other hand, the triangle inequality and Equation
\eqref{eqn:bonneborne} implies that
$$
d_{\infty}(\xi',k_2)\ge d_{\infty}(k_2,\xi)-d_{\infty}(\xi',\xi)\ge
r_i(H)-u=\frac 12-u\;.
$$
Thus we have $r_e(H')< d_{\infty}(\xi',k_2)$ if
\begin{equation}
  2^{-\frac 1a}u^2 < 1/2-u.\label{eqn:c1}
\end{equation}
Similarly, if 
\begin{equation}\label{eqn:c2}
 2^{-\frac 1a}u^2< u-2^{-\frac 1a}e^{-t}
\end{equation}
then the inequality $ r_e(H')< d_{\infty}(\xi',k_1) $
is satisfied.

\medskip
The first condition \eqref{eqn:c1} holds if (and only if)
$$
2^{\frac 1a - 1}\left(-1-\sqrt{1+2^{1-\frac 1a}}\right)
< u <
2^{\frac 1a - 1} \left(-1+\sqrt{1+2^{1-\frac 1a}}\right)
\,.
$$
The lower bound is nonpositive, so this reduces to
\begin{equation}\label{eqn:c1'}
  u< u_1(a)=2^{\frac 1a - 1}
\left(\sqrt{1+2^{1-\frac 1a}}-1\right)\,.
\end{equation}
Let us define $s=2^{-\frac 2a}e^{-t}$.  As $2^{2-\frac 2a}e^{-t_1(a)}<
1$ (by an easy computation), we have $4s<1$. Hence the second
condition \eqref{eqn:c2} holds if (and only if)
\begin{equation}
2^{\frac 1a - 1} \left(1-\sqrt{1-4s}\right) 
< u <
2^{\frac 1a - 1} \left(1+\sqrt{1-4s}\right)
\,.
\label{eqn:2nd}
\end{equation}
Note that $\xi'$ belongs to $\O\,H$, hence by Corollary
\ref{coro:ballshadsize}, we have $u\leq r_e(H)=2^{-\frac 1a}$.  Now,
by an easy computation, the upper bound in \eqref{eqn:2nd} is strictly
bigger than $2^{-\frac 1a}$ if $a<2$, or if $a\geq 2$ and
$e^{-t}<1-2^{-\frac 2a}$. Hence, as $t>t_1(a)$, Equation
\eqref{eqn:2nd} reduces to
\begin{equation}
  u> u_2(a)= 2^{\frac 1a - 1}
\left(1-\sqrt{1-4s}\right)\;.\label{eqn:c2'}
\end{equation}
At least one of the conditions \eqref{eqn:c1'} and \eqref{eqn:c2'}
holds if $u_1(a)> u_2(a)$, that is if
$$
\sqrt{1+2^{1-\frac 1a}}-1>1-\sqrt{1-4s}
\;.
$$
This last equation is equivalent to
$e^{-t}<2^{2/a}\left(\sqrt{1+2^{1-1/a}}-1-2^{-1-1/a}\right)$, which proves
the result in the first case, as $t>t_1(a)$.


\bigskip
{\bf Case 2:} 
Assume that $\xi'$ lies on the left of $k_1$. Then, we assume that $K$
meets $\O\,H'(t)$, and we prove that every point $\eta$ in the right
component of $\overline{\O\,H'-\O\,H'(t)}$ belongs to $K$. By
connectedness, we 
have that $\eta$ lies on the right of $k_1$,
and it is sufficient to prove that $\eta$ lies in $\O\,H$. By
Corollary \ref{coro:ballshadsize}, we only have to prove that
$d_\infty(\xi,\eta)\leq r_i(H)=\frac{1}{2}$.

\medskip
As $H$ and $H'$ bound disjoint open horoballs, as $r_i(H)\leq r_i(H')$
and as $K$
meets $\O\,H'(t)$, the point $\xi'$ belongs to $[\xi,k_1]_\infty$.
Hence by connectedness, the point $\xi'$ belongs to $\O\,H(t)$. By
Corollary \ref{coro:ballshadsize}, we have
$$d_\infty(\xi,\xi')\leq r_e(H(t))=2^{-\frac{1}{a}}e^{-t}\;,$$
and
$$d_\infty(\xi',\eta)\leq r_e(H')\;.$$
By Corollary \ref{coro:disjointhamenstadt}, we have
$r_i(H')=2r_i(H')r_i(H)\leq \frac{1}{2}d_\infty(\xi,\xi')^2$,
so that
$$d_\infty(\xi',\eta)\leq r_e(H')=2^{1-\frac
  1a}r_i(H')\leq 2^{-\frac
  1a}d_\infty(\xi,\xi')^2\leq 2^{-\frac{3}{a}}e^{-2t}\;.$$
By the
triangular inequality, we have
$$d_\infty(\xi,\eta)\leq d_\infty(\xi,\xi')+d_\infty(\xi',\eta)\leq
2^{-\frac{1}{a}}e^{-t}+2^{-\frac{3}{a}}e^{-2t}\;.$$

As $t>t_1(a)$, we have $2^{-\frac{1}{a}}e^{-t}+2^{-\frac{3}{a}}e^{-2t}
\leq\frac{1}{2}$, by an easy computation. Hence $d_\infty(\xi,\eta)\leq
\frac{1}{2}$, which proves the result in the second case.

\bigskip
{\bf Case 3:} 
Finally, assume that $\xi'$ lies to the right of $k_2$.  Then, we
assume that $K$ meets $\O\,H'(t)$, and we prove that the left
component of $\overline{\O\,H'-\O\,H'(t)}$ belongs to $K$. 

\medskip
We start by proving two lemmas, which are obvious in the
constant curvature case, but not immediate otherwise.

\begin{lemma}\label{lem:sensdessusdessous}
  If $H,H'$ are horospheres in $M$ centered at points $\xi,\xi'$
  different from $\xi_\sharp$, bounding disjoint open horoballs, and
  such that there exists a geodesic line $\ell$ starting from
  $\xi_\sharp$ that first meets $H'$ and then meets $H$, then
$$
  r_i(H)\leq 2^{1-\frac{1}{a}}r_i(H')\;.
$$
\end{lemma}

Note that the same result with the constant $2^{1-\frac{1}{a}}$
replaced by $1$ is not true, as shows an easy example built by
deforming a constant curvature $-2$ situation with two tangent
horospheres of the same height.

\medskip \dem Up to pushing $\ell$, we may assume that $\ell$ is
tangent to $H$, and that the point at infinity of $\ell$ different
from $\xi_\sharp$ belongs to $[\xi,\xi']_{\infty}$. Up to replacing
$H'$ by a horosphere also centered at $\xi'$ and contained in the open
horoball bounded by $H'$ (which does not increase $r_i(H')$), we may
assume that $\ell$ is also tangent to $H'$. Orient $\partial
M-\{\xi_\sharp\}$ such that $H$ is on the left of $\ell$. Note that
$H'$ is on the right of $\ell$. If the point
$x$ on $\ell$ moves towards $\xi_\sharp$ on $\ell$, then, with $H_x$
the horosphere tangent at $x$ to $\ell$ and on the left of $\ell$, the
quantity $d(H_x,H_\sharp)$ is non increasing (we assume, as we may,
that $H_\sharp$ is close to $\xi_\sharp$).  Hence, up to replacing $H$
by the horosphere on the left of $\ell$ tangent to $\ell$ at the same
point than $H'$, we may assume that $H$ and $H'$ are tangent. Let $x$
be the (only) tangency point.  Consider the quantity
$$\Delta_H=d(x,H_\sharp)-d(H,H_\sharp)\;.$$
If the curvature of $M$ was constant, equal to $-b^2$ (with $1\leq
b\leq a$), then, as an easy computation in the upper halfspace
model proves, we would have $\Delta_H=\log 2^{\frac{1}{b}}$. Since the
curvature of $M$ lies between $-a^2$ and $-1$, as an easy comparison
argument shows, we have
$$\log 2^{\frac{1}{a}}\leq \Delta_H\leq \log 2\;.$$
By symmetry, this is also true for $H'$. Hence
$r_i(H)/r_i(H')=e^{\Delta_{H'}-\Delta_H}\leq 2^{1-\frac{1}{a}}$,
which proves the result.
\cqfd

\begin{lemma}\label{lem:shadowsinorder} 
  Let $H,H',H''$ be horospheres in $M$ centered at points at infinity
  respectively $\xi,\xi',\xi''$ distinct from $\xi_\sharp$ and with
  $\xi''$ between $\xi$ and $\xi'$ on $\partial M-\{\xi_\sharp\}$,
  such that $H$ is tangent to $H''$ and the open horoball bounded by
  $H$ is disjoint from both open horoballs bounded by $H'$ and $H''$.
  Assume that there exists no geodesic line  starting from
  $\xi_\sharp$ that first meets $H'$ and then meets $H$. Then
$$
  \O\,H'\cap[\xi,\xi'']_\infty\subset \O\,H''\;.
$$
\end{lemma}

\dem
Assume by absurd that there exists a point in
$(\O\,H'-\O\,H'')\cap[\xi,\xi'']_\infty$. Let $\ell$ be the geodesic
line starting from $\xi_\sharp$, which is tangent to $H''$, and whose
point at infinity lies in $[\xi,\xi'']_\infty$, which exists as $H\cup
H''$ separates $M$.  By convexity of horoballs, the complement in
$\ell$ of the open horoball bounded by $H$ is the union of two
geodesic rays $\ell_1,\ell_2$ with $\ell_2$ ending in
$[\xi,\xi'']_\infty$.  Since $\ell$ cannot first meet $H'$ and then meet
$H$, as $H,H'$ are bounding disjoint open horoballs, and by convexity,
the horosphere $H'$ meets $\ell_2$ in two distinct points $x,y$.

\begin{figure}[ht]
\input{newimpossible.pstex_t}
\caption{}\label{fig:impossible}
\end{figure}

Assume that $x$ is closer to $\xi_\sharp$ on $\ell$, and orient the
topological line $H'$ from $y$ to $x$. Let $p$ be the first
intersection point of $H'$ with $H''$ after $x$, and $\rho'',\rho'$
the geodesic rays starting from $p$ and ending at $\xi'',\xi'$.  As
$\xi'$ is on the right of $\xi''$ on $\partial M-\{\xi_\sharp\}$, and
as $H'$ enters the horoball bounded by $H''$ at $p$, the tangent
vector $v'$ to $\rho'$ at the origin is contained in the angular
sector between the tangent vector of $H'$ and the tangent vector $v''$
to $\rho''$ at the origin.  This contradicts the fact that an
horosphere is perpendicular to the geodesic lines ending at the point
at infinity of a horosphere.  
\cqfd

\medskip Let us end now the proof of Case 3. Consider the (unique by
convexity) horosphere $H''$ centered at $k_2$ and tangent to $H$.  Let
$\ell_{k_2}$ be the geodesic line starting from $\xi_\sharp$ tangent
to $H$ at a point $x$ and ending at $k_2$.
As $K$ meets $\O\,H'(t)$, the horoball bounded by $H'(t)$ has to meet
$\ell_{k_2}$ in a segment. If $\ell_{k_2}$ was first meeting $H'(t)$
and then $H$, then by Lemma \ref{lem:sensdessusdessous}, we would have
$$
r_i(H)\leq 2^{1-\frac{1}{a}}r_i(H'(t))\;.
$$
As $0<r_i(H')\leq r_i(H)$, we would have $e^{-t}\geq
2^{\frac{1}{a}-1}$. This is impossible, as $t> t_1(a)$.

Hence no geodesic line starting from $\xi_\sharp$ first meets $H'$ and
then meets $H$. By Lemma \ref{lem:shadowsinorder}, we have
$\O\,H'\cap[\xi,k_2]_\infty\subset \O\,H''\cap[\xi,k_2]_\infty$.
Hence the result follows from Case 1 (applied by replacing $H'$ by
$H''$.)
\cqfd

\bigskip
\noindent{\bf Proof of Theorem \ref{thm:2dimbounded}. } 
Theorem \ref{thm:2dimbounded} follows from Proposition
\ref{prop:2dimhoroballs} in the same way as Theorem \ref{thm:bounded2}
follows from Proposition \ref{prop:twoabstshad}.  Note that the
assumption about the space being geodesic and having extendable
spheres was only used in Proposition \ref{prop:twoabstshad}, and not
in the deduction of Theorem \ref{thm:bounded2} from Proposition
\ref{prop:twoabstshad}. Note that in the proof, we fix a connected
component $K_0$ of $\overline{\O\,H_0-\O\,H_0(t)}$ and we find a
geodesic line starting from $\xi_\sharp$ as wanted, which furthermore
has its endpoint in $K_0$.  But there exists another such connected
components $K'_0$ (which is different from $K_0$) and working with
$K'_0$ instead of $K_0$, we get a second geodesic line as wanted.
\cqfd

\medskip
\begin{remark}\label{rem:2-dim}
  (a) Theorem \ref{thm:2dimbounded} gives the value
  $t_{1}(1)=-\log(4\sqrt 2-5)\approx 0.42$ for the case of constant
  curvature $-1$, improving the value $t_0(1)=-\log(\sqrt 5-2)\approx
  1.44$ given by Theorem \ref{thm:bounded2}.  
%

\medskip
\noindent(b) If $M=\hdr$, we have seen in the proof that we may
replace the condition $t>t_1$ by $t\geq t_1$, and the triangle
inequalities used in the proof of Proposition \ref{prop:2dimhoroballs}
are replaced by equalities. In fact Theorem \ref{thm:2dimbounded} is
sharp, in the following sense: there exists a family of horoballs
$(H_n)_{n\in\NN}$ in $\hdr$ (with the upper halfspace model) such that
there exists a geodesic line starting at $\infty$ and meeting the open
horoball bounded by $H_0$ which avoids the open horoballs bounded by
$H_n(t)$ for every $n$ in $\NN$ if and only if $t\geq
t_{1}(1)=-\log(4\sqrt 2-5)$. To construct such a family, we start with
a horosphere $H_0$ (not centered at $\infty$).  We define $H_1, H_2$
to be the boundaries of the two maximal horoballs whose interiors are
disjoint from the open horoball bounded by $H_0$, and which are
contained in the shadow of $H_0$ seen from $\infty$. An easy
computation shows that $\O\,H_1$ and $\O\,H_2$ are exactly the two
components of $\overline{\O\,H-\O\,H(t_1(1))}$. And then one iterates
the construction (see picture below).
 
\begin{figure}[ht]
\input{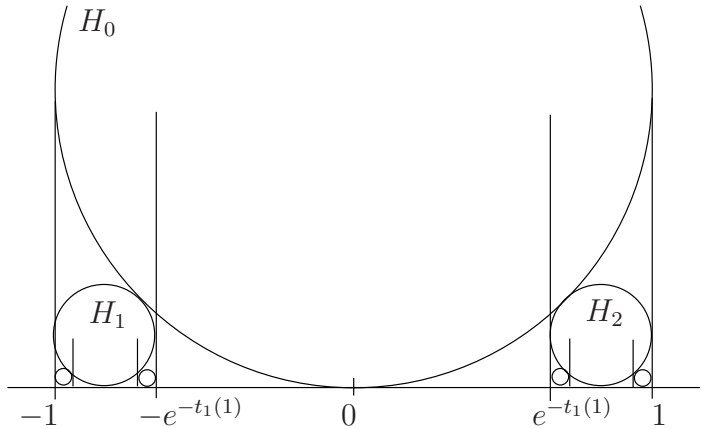}
\caption{An extremal packing for Theorem
\ref{thm:2dimbounded}.}\label{fig:extremal}
\end{figure}

\medskip
\noindent(c) Let $(H\!B_r)_{r\in\QQ\cup\{\infty\}}$ be the 
Apollonian packing in the upper half plane model of $\hdr$ as in
Remark \ref{rem:postabstract}(ii). 
Let $\xi\in\RR$. If $[\xi,\infty]$ meets $H\!B_{p/q}(t)$, then
\begin{equation}
  \left|\xi-\frac pq\right|<\frac{e^{-t}}{2 q^2}. \label{eqn:dirichlet}
\end{equation}
It is well known (see 
\cite{Ford, HW, Khi}) 
that there are irrational numbers $\xi$ such that Equation
\eqref{eqn:dirichlet} has only finitely many solutions $\frac pq$ if
$t>-\log\frac 2{\sqrt 5}\approx 0.11$. Furthermore, if
$\xi=\frac{1+\sqrt{5}}{2}$ is the golden ratio, then by \cite{RWT},
Equation \eqref{eqn:dirichlet} has no solution at all if $t>-\log(3-{\sqrt
  5})\approx 0.269$. Theorem \ref{thm:2dimbounded} implies that there
exists at least one irrational number $\xi$ such that
\eqref{eqn:dirichlet} has no solution for $t\ge t_1(1)=-\log(4\sqrt
2-5)\approx 0.42$.  This is not much more than what is actually gotten
for the golden ratio.

\medskip
\noindent(d) Clearly, the two possible choices of the maximal
interval, in the analog of the first step of the proof of Theorem
\ref{thm:abstract} within the proof of Theorem \ref{thm:2dimbounded},
give two different geodesics which avoid the scaled horospheres such
that the endpoints of these geodesics are both in the shadow of the
first horosphere used in the construction.
\end{remark}

\section{Unclouding the sky of real hyperbolic spaces}
\label{sec:hnr}

Let $\hnr$ be the real hyperbolic $n$-space (with constant curvature
$-1$), for $n\geq 2$. In this section, we will use the symmetries of
$\partial\hnr$ to show that the condition $t\geq t_1(1)= -\log
(4\sqrt{2}-5)$, which is our best estimate for the real hyperbolic
plane (see Theorem \ref{thm:2dimbounded} and Remark \ref{rem:2-dim}
(a)), also works in the higher-dimensional situation.

\begin{theorem}\label{thm:hnrbounded} 
  Let $\xi_\sharp$ be a point in $\partial \hnr$, and $H_\sharp$ be a
  horosphere in $\hnr$ centered at $\xi_\sharp$.  For every sequence
  $(H_n)_{n\in\NN}$ of horospheres in $\hnr$ bounding open horoballs
  that are pairwise disjoint and disjoint from the open horoball
  bounded by $H_\sharp$, if $t\geq t_1(1)$, then there exists at least
  two geodesic lines in $\hnr$, starting at $\xi_\sharp$, which avoid
  $H_n(t)$ for every $n$ in $\NN$.
\end{theorem}

\dem By homogeneity, we may use the upper halfspace model for $\hnr$,
with $\xi_\sharp$ the point at infinity $\infty$ and $H_\sharp$ the
horizontal Euclidean hyperplane at Euclidean height $1$. In
particular, as we have already seen, the space
$\partial\hnr-\{\infty\}$ is the horizontal coordinate hyperplane
$\RR^{n-1}$, the shadows of horoballs (seen from $\infty$) are
Euclidean balls in $\RR^{n-1}$, and the Hamenst\"adt distance is the
Euclidean distance.

\begin{proposition}\label{prop:hnrhoroballs}  
  Let $H$ and $H'$ be horospheres in $\hnr$ bounding open horoballs
  that are disjoint and disjoint from the open horoball bounded by
  $H_\sharp$, with $H'$ lower than $H$ (i.e.~$r_i(H')\leq r_i(H)$).
  Let $t\geq t_{1}(1)$.  For every maximal Euclidean ball $K$ of
  $\overline{\O\,H-\O\,H(t)}$, the following dichotomy holds~:
\begin{itemize}
\item either $K$ does not meet $\O\,H'(t)$
\item or there exists a maximal Euclidean ball $K'$ in
  $\overline{\O\,H'-\O\,H'(t)}$, which is contained in $K$.
\end{itemize}
\end{proposition}

\dem Assume that $K$ meets $\O\,H'(t)$.
Let $x$ be the common center of the Euclidean balls $\O\,H(t)$, let $x'$
be the common center of the Euclidean balls $\O\,H'(t)$, let $y$ be
the center of $K$ (note that $y\neq x$), let $D$ be the Euclidean line
passing through $x$ and $y$, and let $D'$ be the Euclidean line
through $x'$ and $y$ (take $D'=D$ if $x'=y$).

\medskip
\noindent
{\bf Case 1: } Assume first that $D=D'$ (see the picture below).


\begin{figure}[ht]
\input{fig_jouni4jan03.pstex_t}
\caption{}\label{fig__jouni4jan03}
\end{figure}

Let $A$ be the hyperbolic plane in $\hnr$ with boundary $D\cup
\{\infty\}$. Note that $H\cap A$ and $H'\cap A$ are horoballs in $H$,
and that $\O\,H\cap D$ is the shadow of $H\cap A$ seen from $\infty$
in $A$. Note that, by maximality, $K\cap D$ is a connected component
of $\left(\overline{\O\,H-\O\,H(t)}\right)\cap D$. Also note that
$K\cap D$ meets $\O\,H'(t)\cap D$. Hence, by Proposition
\ref{prop:2dimhoroballs} applied to $M=A$, there exists a connected
component $I$ of $\left(\overline{\O\,H'-\O\,H'(t)}\right)\cap D$
which is contained in $K\cap D$.  Consider the ball $K'$ in $\partial
\hnr-\{\infty\}$ having $I$ as one diameter. Then $K'$ is contained in
$\overline{\O\,H'-\O\,H'(t)}$ and is maximal there, since $D$ goes
through the center of the balls $\O\,H'(s)$. Furthermore $K'$ is
contained in $K$, since the ball with diameter a segment contained in
$K$ of a line passing through the center of $K$ is contained in $K$.

\medskip
\noindent {\bf Case 2: } 
Now, assume that $D'$ and $D$ are different (see the picture below).

Let $P$ be the orthogonal subspace of $\RR^{n-1}$ through $y$ to the
plane containing $D$ and $D'$. Let $f$ be the isometry of $X$ fixing
$\infty$, preserving the horosphere $H_\sharp$, and inducing on
$\partial\hnr-\{\infty\}=\RR^{n-1}$ the Euclidean rotation fixing $P$
and sending $D'$ to $D$ such that $y$ is between $x$ and $x'$.  Then
the map $f$ preserves $K$, and sends $H'$ to a horoball $H''$ which
is still disjoint from (and lower than) $H$. By the case $D=D'$, let
$K''$ be a maximal Euclidean ball in $\overline{\O\,H''-\O\,H''(t)}$,
which is also contained in $K$. Then, since $f$ is an isometry of $X$,
the Euclidean ball $f^{-1}(K'')$ is a maximal Euclidean ball in
$\overline{\O\,H'-\O\,H'(t)}$, which is also contained in
$f^{-1}(K)=K$.
\cqfd
\medskip

\begin{figure}[ht]
\input{fig_jouni2_4jan03.pstex_t}
\caption{}\label{fig_jouni2_4jan03}
\end{figure}

\medskip
Now, Theorem \ref{thm:hnrbounded} follows from the above Proposition
\ref{prop:hnrhoroballs} exactly in the same way as Theorem
\ref{thm:2dimbounded} followed from Proposition
\ref{prop:2dimhoroballs}.  \cqfd

\section{Bi-infinite geodesics, geodesic rays, and closed geodesics 
in finite volume cusped manifolds}
\label{sec:bi-infinite}

In this concluding section, we give some applications of the results
obtained in the previous sections. We start by studying the problem of
finding geodesic lines which avoid the scaled family of horoballs in
both directions, as well as geodesic rays starting from given points.

\begin{theorem}
  Let $X$ be a proper geodesic CAT($-1$) space, and $t_{\rm
  min}>0$. Assume that one of the following conditions holds:
\begin{enumerate}
\item $X$ is the real hyperbolic $n$-space with $n\geq 2$, and
  $t_{\rm min}=t_1(1)=-\log(4\sqrt{2}-5)$;
\item $X$ is a complete simply connected Riemannian manifold of
  dimension $2$, with pinched 
  curvature $-a^2\leq K\leq -1$,
  and $t_{\rm min}=t_1(a)$;
\item $X$ is a locally finite metric tree, without vertices of degree
  $1$ or $2$, with edge lengths at most $\ell_{\rm max}$, and $t_{\rm
    min}=\ell_{\rm max}$;
\item $X$ is a  symmetric space, and $t_{\rm min}=t_0$ (see
Theorem \ref{thm:bounded2}).
\end{enumerate}
Let $(H_n)_{n\in\NN}$ be a sequence of horospheres in $X$ bounding
pairwise disjoint open horoballs $(H\!B_n)_{n\in\NN}$. If $t> t_{\rm
  min}$ in case (3) and $t> \log(2+\sqrt{5})+ t_{\rm min}$ otherwise,
for every $x$ in $X-\bigcup_n \stackrel{\circ}{H\!B_n}$, then there
exists at least one geodesic ray starting at $x$ which avoids $H_n(t)$
for every $n$ in $\NN$.

Furthermore, if $t>t_{\rm min}$ in Case (3) and $t> t_{\rm min}
+\log(1+\sqrt{2})$ otherwise, then there exists at least one geodesic
line which completely avoids $H_n(t)$ for every $n$ in $\NN$.
\end{theorem}

The constants in this theorem are probably not optimal (except for the
case (3), where they are optimal). The first assertion of this theorem
implies the first assertion in Theorem \ref{theointro} in the
introduction. The second assertion of Theorem \ref{theointro} follows
from Theorem \ref{thm:hnrbounded}, Theorem \ref{thm:2dimbounded},
Remark \ref{rem:finquatre} and Theorem \ref{thm:bounded2}.

\medskip
\dem Let $(H_n)_{n\in\NN}$ be as in the statement, and define $\xi_n$
to be the point at infinity of $H_n$.

\medskip
Let us prove the first assertion. Let $x,t$ be as in the statement.
Let $H_{n_0}$ be one of the closest horoballs to $x$ (which exists
since the family of horoballs is locally finite).  Let $H_{n_1}$ be
the first horosphere met after $x$ by the geodesic line starting from
$\xi_0$ and passing through $x$, with $y$ the intersection point. If
$H_{n_1}$ does not exists, then we are done. We may assume that $n_0=0$
and $n_1=1$.

\blemm 
For every point at infinity $\xi_\sharp$ and every horosphere
$H$ not centered at $\xi_\sharp$ in a CAT$(-1)$ space, define
$\C^-_{\xi_\sharp}H$ to be the union of the geodesic rays meeting $H$
exactly in their starting point, and converging to
$\xi_\sharp$. For every such ray $\rho$, the subset
$\C^-_{\xi_\sharp}H$ is contained in the
$\log(2+\sqrt{5})$-neighborhood of $\rho$. 
\elemm

\dem By comparison, this follows from an easy computation in the 
upper halfspace model of the real hyperbolic plane. Note that
in a CAT$(-a^2)$ space (resp.~a tree), we may replace the constant 
$\log(2+\sqrt{5})$ by $\frac{1}{a}\log(2+\sqrt{5})$ (resp.~$0$).
\cqfd

\medskip Define $\alpha=0$ in the case (3) and
$\alpha=\log(2+\sqrt{5})$ otherwise. By construction, no horosphere
$H_n$ besides $H_0$ meets the geodesic ray from $y$ to $\xi_0$. Hence
by the lemma, no horosphere $H_n(s)$ with $n\geq 1$ and $s> \alpha$
meets $\C^-_{\xi_0}H_1$.  As $t-\alpha>t_{\rm min}$, we can apply
Theorem \ref{thm:hnrbounded}, Theorem \ref{thm:2dimbounded}, Remark
\ref{rem:finquatre} and Theorem \ref{thm:bounded2} respectively, while
taking $H_1$ as the starting horosphere in our inductive
constructions. Hence there exists at least one geodesic line $\ell$
starting from $\xi_0$ which avoids $H_n(t-\alpha)$ for $n\geq 1$, and
whose other endpoint belongs to $\O_{\xi_0}H_1$. Note that in Case
(4), we may indeed take $\xi_\sharp=\xi_0$ to apply Theorem
\ref{thm:bounded2}, as the isometry group of $X$ then acts
transitively on $\partial X$.

By the above lemma, the point $x$ lies at distance at most
$\alpha$ from the geodesic line $\ell$. Hence by convexity,
the geodesic ray $\rho_0$ starting from $x$ and converging to the
point at infinity of $\ell$ different from $\xi_0$ is contained in the
$\alpha$-neighborhood of $\ell$. Hence $\rho_0$ avoids
$H_n(t)$ for $n\geq 1$. The result follows, because, by construction,
$\rho_0$ also avoids $H_0$.

\medskip Now, let us prove the second assertion.  Let $t$ be as in the
statement.  By Theorem \ref{thm:hnrbounded}, Theorem
\ref{thm:2dimbounded}, Remark \ref{rem:finquatre} and Remark
\ref{rem:doublerlamise} respectively, there exist at least two
geodesic lines starting from $\xi_0$ which, for $n\geq 1$, avoid
$H_n(t)$ in Case (3) and $H_n(t-\log(1+\sqrt{2}))$ otherwise, and are
contained in $\O_{\xi_0}H_{n_0}$ for some $n_0$ by construction.  Note
that in Case (4), we may indeed take $\xi_\sharp=\xi_0$ to apply
Remark \ref{rem:doublerlamise}, as the isometry group of $X$ then acts
transitively on $\partial X$, and as the spheres in $\partial
X-\{\xi_0\}$ have antipodal points (see Remark \ref{rem:CCantipod}).

Let $\eta,\eta'$ be the endpoints of these two geodesic lines. By
Lemma \ref{lem:lowgeod}, the geodesic line between $\eta$ and $\eta'$
does not intersect the open horoball bounded by $H_0$, and in
particular avoids $H_0(t)$.

In case (3), the geodesic line between $\eta$ and $\eta'$ is contained
in the union of the geodesic lines between $\xi_0$ and $\eta$ and
between $\xi_0$ and $\eta'$, which proves the result.

For the other cases, recall that in a CAT$(-1)$ space, given three
points at infinity, and three geodesic lines between the pairs of
them, each geodesic line is contained in the
$\log(1+\sqrt{2})$-neighborhood of the union of the two others. This
follows, by comparison with an ideal triangle in the hyperbolic plane,
from an easy computation in the upper halfplane model (see also
\cite[Proposition 2.21]{GH}, where the worse constant $\log 3$ is
given). Hence the geodesic line between $\eta$ and $\eta'$ also avoids
$H_n(t)$ for $n\geq 1$.
\cqfd

\brema {\rm Let $X$ be a proper geodesic CAT($-1$) space, $\xi_\sharp$
  a point at infinity, $H_\sharp$ a horosphere in $X$ centered at
  $\xi_\sharp$, and $t_{\rm min}>0$. With $d_\ell$ the length distance
  associated to the Hamenst\"adt distance $d_{\xi_\sharp,H_\sharp}$,
  assume that there exists $C\geq 1$ such that
  $d_{\xi_\sharp,H_\sharp}\leq d_\ell\leq Cd_{\xi_\sharp,H_\sharp}$,
  the metric space $(\partial X-\{\xi_\sharp\},d_\ell)$ has extendable
  spheres with modulus $\delta(\cdot)$, and the spheres in $(\partial
  X-\{\xi_\sharp\},d_\ell)$ have antipodal points.  Then for every
  sequence $(H_n)_{n\in\NN}$ of horospheres bounding open horoballs
  that are pairwise disjoint and disjoint from the open horosphere
  bounded by $H_\sharp$, if $t> t_0(C,\delta) +\log(1+\sqrt{2})$, then
  there exists at least one geodesic line which completely avoids
  $H_\sharp$ and $H_n(t)$ for every $n$ in $\NN$.
} 
\erema

\dem 
The proof is the same as the previous one, we just added the
hypothesis necessary to handle the more general case. \cqfd

\bigskip 
Let us now give an application of our results to equivariant
families of horospheres.
Let $V$ be a complete 
 noncompact nonelementary geometrically finite Riemannian manifold
with sectional curvature $K\le -1$.
Let $\beta_{e}$ be the Busemann function on $V$ with respect to a cusp
$e$ of $V$, normalised to be $0$ on the boundary of the maximal open
Margulis neighbourhood $N_e$ of $e$ in the convex core of $V$.  Recall
that (see for instance 
\cite{HP3}) if $\rho:[0,+\infty)\ra V$ is a geodesic ray contained in
the closure of $N_e$ with $\rho(0)$ not in $N_e$, then
$\beta_e(x)=\lim_{t\ra\infty}(t-d(\rho(t),x))$. The height with
respect to $e$ of any compact subset $A$ in $V$ can be defined as
$$
\height_e(A)=\max_{x\in A}\beta_{e}(x).
$$
We denote by $h_e(V)$ the infimum of the heights of the closed
geodesics on $V$ with respect to $e$.  The results of the previous
sections give an upper bound on $h_e(V)$.

\btheo\label{thm:finite vol}
  Let $V$ be a complete noncompact
nonelementary geometrically finite 
Riemannian manifold (for instance a manifold with finite volume), 
   $e$ be a cusp of $V$, and $t_{\rm min}>0$. 
 Assume that one of the following conditions holds:
\begin{enumerate}
\item $V$ has constant sectional curvature $-1$, and
  $t_{\rm min}=t_1(1)=-\log(4\sqrt{2}-5)$;
\item $V$ is two-dimensional with pinched 
  curvature $-a^2\leq K\leq -1$, and
  $t_{\rm min}=t_1(a)$;
\item  the universal cover of $V$ is a CAT($-1$) space and
satisfies the conditions (i) and (ii) of Section \ref{sec:spaboundsha}
with $\xi_\sharp$ a lift of the cusp $e$, and $t_{\rm min}=t_0$
(see Theorem \ref{thm:bounded2}).
\end{enumerate}
Then $h_e(V)\leq t_{\rm min}$.
\etheo

This theorem
implies Corollary \ref{corointro} in the introduction.
Note that if a negatively curved homogeneous Riemannian manifold has a
finite volume quotient, then it is symmetric \cite{Hei2}. Hence, in
Corollary \ref{corointro}, we
only gave the result for locally symmetric spaces.
For related results on bounded geodesics in Riemannian manifolds, see
\cite{Sch} and the references therein.

\medskip
\dem 
The lift, to the universal cover $\wt V$ of $V$, of the union of the
maximal open Margulis neighbourhoods of ends of $V$ is a disjoint
collection of open horoballs in $\wt V$. Thus, the theorems
\ref{thm:hnrbounded}, \ref{thm:2dimbounded}, \ref{thm:bounded2}
respectively can be applied to the corresponding family of
horospheres.  By \cite[Theorem 3.4]{HP2}, $h_e(V)$ is equal to the
lower bound of all $h$ in $\RR$ such that there exists a geodesic line
starting from $e$, which does not converge into a cusp of $V$, and
eventually avoids $\beta_e^{-1}([h,+\infty))$. This implies the
result.  
\cqfd

\bigskip
{\small\noindent \begin{tabular}{l} 
Department of Mathematics and Statistics\\
 P.O. Box 35\\ 
40014 University of Jyv\"askyl\"a, FINLAND.\\
{\it e-mail: parkkone@maths.jyu.fi}
\end{tabular}
\medskip

\noindent \begin{tabular}{l}
D\'epartement de Math\'ematique et Applications, UMR 8553 CNRS\\
Ecole Normale Sup\'erieure\\
45 rue d'Ulm\\
75230 PARIS Cedex 05, FRANCE\\
{\it e-mail: Frederic.Paulin@ens.fr}
\end{tabular}
}

\end{document}